\documentclass[11pt]{article}
\usepackage[utf8]{inputenc}
\usepackage[margin=1in,left=1in]{geometry}
\usepackage{setspace}
\usepackage{graphicx}
\usepackage{amssymb, amsmath, amsthm, graphics}
\usepackage[mathscr]{euscript}
\usepackage{tikz}
\usepackage{pigpen}
\usepackage{bm}
\usetikzlibrary{calc}
\usetikzlibrary{matrix}
\usepackage{mathtools} 

\usepackage{enumitem}

\usepackage{latexsym}
\usepackage[all]{xy}
\usepackage{color}

\usepackage{bbold}

\usepackage{tikz}
\input xy
\xyoption{all}
\usepackage{euscript}
\usepackage{multirow}
\usepackage{etex, pictexwd,dcpic}
\usetikzlibrary{positioning}
\usetikzlibrary{shapes.geometric}
\usetikzlibrary{shapes.misc}
\usetikzlibrary{calc}
\usetikzlibrary{positioning}

\usepackage{pgflibraryarrows}
\usepackage{pgflibrarysnakes}

\usetikzlibrary{trees} 
\usetikzlibrary[trees] 

\newcommand{\divides}{\bigm|}
\newcommand{\ndivides}{%
  \mathrel{\mkern.5mu 
    \ooalign{\hidewidth$\big|$\hidewidth\cr$\nmid$\cr}%
  }%
}

\newtheorem{theorem}{Theorem}
\newtheorem{definition}[theorem]{Definition}
\newtheorem{lemma}[theorem]{Lemma}
\newtheorem{corollary}[theorem]{Corollary}
\newtheorem{proposition}[theorem]{Proposition}

\newtheorem{remark}{Remark}

\newtheorem{example}{Example}

\numberwithin{equation}{section}

\renewcommand{\(}{\begin{equation*}}
\renewcommand{\)}{\end{equation*}}
\newcommand{\bea}{\begin{eqnarray*}}
\newcommand{\eea}{\end{eqnarray*}}

\def\endofproof {\hfill{$\Box$}\\}

\newcommand{\cA}{\ensuremath{\mathcal A}}

\newcommand{\beq}{\begin{equation}}
\newcommand{\eeq}{\end{equation}}

\newcommand{\into}{\hookrightarrow}

\newcommand{\theproof}{\noindent {\bf Proof.\ }}

\numberwithin{equation}{section}

\renewcommand{\(}{\begin{equation}}
\renewcommand{\)}{\end{equation}}

\newcommand{\RP}{\RR \text{P}}

\def\ch{{\rm  ch}}

\def\1{{\bf 1}}

\def\<{\langle}
\def\>{\rangle}

\numberwithin{equation}{section}


\newcommand{\R}{\ensuremath{\mathbb R}}
\newcommand{\RR}{\ensuremath{\mathbb R}}
\newcommand{\NN}{\ensuremath{\mathbb N}}

\newcommand{\ZZ}{\ensuremath{\mathbb Z}}
\newcommand{\Z}{\ensuremath{\mathbb Z}} 
\newcommand{\QQ}{\ensuremath{\mathbb Q}}

\newcommand{\BB}{\ensuremath{\mathbf B}}

\newcommand{\chp}{\ensuremath{\mathscr{C}\mathrm{h}^{+}}}

\newcommand{\sset}{\ensuremath{s\mathscr{S}\mathrm{et}}}

\newcommand{\set}{\ensuremath{\mathscr{S}\mathrm{et}}}

\newcommand{\sab}{\ensuremath{s\mathscr{A}\mathrm{b}}}
\newcommand{\ab}{\ensuremath{\mathscr{A}\mathrm{b}}}

\newcommand{\cartsp}{\mathscr{C}\mathrm{art}\mathscr{S}\mathrm{p}}

\newcommand{\D}{\ensuremath{\mathscr{D}}}

\newcommand{\map}{\mathrm{Map}}

\begin{document}

\title{Primary operations in differential cohomology}

 \author{
 Daniel Grady and Hisham Sati\\
  }

\maketitle

\begin{abstract} 
We characterize primary operations in differential cohomology via stacks,
and illustrate by differentially refining Steenrod squares and Steenrod powers  explicitly. 
This requires a delicate interplay between integral, rational, and mod $p$ cohomology, 
as well as cohomology with $U(1)$ coefficients and differential forms. 
Along the way we develop computational techniques in differential cohomology, including 
a K\"unneth decomposition, that should also be useful in their own right, and point to
applications to higher geometry and mathematical physics. 
 \end{abstract}

 \tableofcontents
 
\section{Introduction}

Cohomology operations with coefficients in a group $G$ are natural transformations
of the form $H^n(-; G) \to H^m(-; G)$. By Brown representability and the Yoneda 
lemma this is equivalent to calculating the universal cohomology group of Eilenberg-MacLane spaces 
$H^m(K(G,n); G)$, which is in turn equivalent 
to calculating the homotopy classes of maps $[K(G, n), K(G, m)]$. Thus all cohomology
operations of fixed degree are accounted for by calculating the cohomology of the 
Eilenberg-MacLane space $K(G, n)$. To account for all cohomology operations
one obviously has to vary both $m$ and $n$. See \cite{MT} 
\cite{St} for detailed accounts. 

\medskip 
We are interested in differential cohomology (see \cite{CS} \cite{Fr} \cite{HS} \cite{SSu} 
\cite{Bun1} \cite{BS} \cite{Urs} \cite{BB}).
What replaces Eilenberg-MacLane 
spaces are various stacks of higher $U(1)$-bundles ($n$-bundles) with connections. Thus, cohomology
operations will involve the differential cohomology of such stacks, and this 
process can be described via mapping spaces of stacks. 
For differential refinements we will need 
to study morphisms of stacks
\(\label{381}
\widehat{\theta}:\BB^nU(1)_{\nabla}\to \BB^mU(1)_{\nabla}\;,
\)
where $\BB^nU(1)_{\nabla}$ represents the moduli stack of $n$-bundles equipped with 
connection, studied in \cite{FSSt}\cite{FSS1}\cite{FSS2}. The homotopy classes of such morphisms will in turn
be describe the differential cohomology group 
$$
\widehat{H}^{k+1}(\BB^{n}U(1)_{\nabla};\ZZ)
:=\pi_0\map(\BB^{n}U(1)_{\nabla},\BB^{k}U(1)_{\nabla})\;.
$$
One of the main goals of this paper is to characterize 
 this group for various values of $k$ and $n$. This in turn will lead to 
 to a full
characterization of primary cohomology operations in differential 
cohomology.

\medskip
Since differential cohomology operations, as we will see, 
 involve various coefficients, we find it useful to point out the interrelations
 that already exist between these (we found the discussion in 
 \cite{FFG} particularly useful). 
 This should also help us develop some intuition 
 for the full differential case.   
 Note that for coefficients being one of $\Z, \Z/p$ or $\QQ$ i.e. an abelian group, then 
the set of all cohomology operations 
$H^m(K(G, n); G)$, where $G$ and $G'$ are from the above 
set, will also be abelian.

 \paragraph{Operations from $\Z/p$ to $\QQ$.}
We know that $H^q(K(G, n); \QQ)=0$ for all $q>0$ when $G$ is a finite abelian group, i.e. for us $\Z/p$. This shows that there are
nontrivial cohomology operations from $\Z/p$-coefficients to 
$\QQ$-coefficients.

 \paragraph{Operations from $\Z$ to $\QQ$.}
 We will distinguish the odd and even cases. For the first, 
 $H^q(K(\Z, 2n+1); \QQ)$ is nonzero only for $q=2n+1$, where it is 
equal to $\QQ$, with generator the image of the fundamental class
$\iota$ under the homomorphism $r: H^{2n+1}(K(\Z, 2n+1); \Z)
\to H^{2n+1}(K(\Z, 2n+1); \QQ)$ induced by the natural embedding
$\Z \hookrightarrow \QQ$. Thus, every operation from an odd-dimensional 
integral class to rational cohomology preserves the dimension, i.e. 
takes $\alpha \in H^{2n+1}(X; \Z)$ to $\lambda r(\alpha) \in H^{2n+1}(X; \QQ)$ for 
some fixed rational number $\lambda$ corresponding to the operation. 
In the even case, $H^q(K(\Z, 2n); \QQ)=\QQ[r(\iota)]$, so that every operation from 
even integral cohomology to rational cohomology is given as 
the power $\alpha \mapsto \lambda \alpha^k$, where 
$k \in \Z$, $\lambda \in \QQ$ are determined by the operation. 

 \paragraph{Operations from $\QQ$ to $\QQ$.}
The group $H^n(K(\QQ, m); \QQ)$ can be straightforwardly calculated via e.g. the Serre spectral sequence. In this case, one has that any cohomology operation assigns to an element $\alpha \in H^n(X; \QQ)$ 
the element $\lambda \alpha^k \in H^{nk}(X; \QQ)$, where
$k \in \Z$, $\lambda \in \QQ$ both fixed by the operation. 

\paragraph{Operations from $H^*(-)_{\rm dR}$ to $H^*(-)_{\rm dR}$.}             
On the other hand, operations in de Rham cohomology can be 
deduced from those on rational cohomology via the de Rham theorem. 
 Hence, de Rham operations
should be systematically characterized. 
The de Rham cohomology groups and the rational cohomology groups
have the same underlying algebraic structure.

 \paragraph{Operations from $\Z/p$ to $\Z/p$.} From an algebraic and homotopic 
 point of view, these are perhaps most studied. 
We start with degree-preserving operations. From $H^n(K(\Z/p; n); \Z/p)=\Z/p$ it follows 
that any degree-preserving such 
operation is multiplication by a scalar in $\Z/p$. Then from 
$H^{n+1}(K(\Z/p, n); \Z/p)=\Z/p$ it follows that there exist a unique operation 
raising the degree by one generating all such operation. This generator is given by the connecting homomorphism $\beta_p$, i.e. the Bockstein homomorphism for the the coefficients 
sequence
$\Z/p \overset{\times p}\longrightarrow \Z/p^2 \overset{\rho_p}\longrightarrow \Z/p$.

\medskip
Note that since $H^q(K(\Z/p, n); \Z/p)=0$ for $n+1 < q < n+2p-2$, there 
are no operations that raise the degree by $2,3,4, \cdots, 2p-3$. 
The next degree where a nonzero operation exits is in dimension 
$2p-1$, where there is a unique operation corresponding to 
$
H^{n+ 2p-2}(K(\Z/p, n); \Z/p)=\Z/p
$
for $n>p$, which is the reduced Steenrod power $P^1$.
The next degrees 
$H^{n+ 2p-1}((K(\Z/p, n); \Z/p)=\Z/p$ and 
$H^{n+2p}(K(\Z/p, n); \Z/p)=\Z/p \oplus \Z/p$ correspond to
combinations of the operations $\beta_p P^1$, $P^1 \beta_p$ and 
$\beta_p P^1 \beta_p$. The next nontrivial degree is $4p-4$
corresponding to the reduced power $P^2$. 

\medskip
The complete classification of these operations for $\Z/p$ is 
given by studying the mod $p$ Steenrod algebra $\cA_p$
(see \cite{St} \cite{C} \cite{Mi} \cite{May}).

\paragraph{Operations from $\Z$ to $\Z/p$.}
The main example here is 
$\rho_p: \Z \to \Z/p$, the mod $p$ reduction, for $p$ a prime number. This
induces an operation of the same name on cohomology
$\rho_p: H^n(-; \Z) \to H^n(-; \Z/p)$.

\paragraph{Operations from $\Z/p$ to $\Z$.}
 Consider $\beta$,  the connecting 
homomorphisms, i.e. the Bockstein homomorphism, for the the coefficients 
sequence
 $\Z \overset{\times p}\longrightarrow \Z \overset{\rho_p}\longrightarrow \Z/p$.
For a class $x \in H^n(X; \Z/p)$, the class $\beta (x)$ is an integral
element of $H^{n+1}(X; \Z/p)$, i.e. it belongs to the image of the 
mod $p$ reduction homomorphism
$\rho_p: H^{n+1}(X; \Z) \to H^{n+1}(X; \Z/p)$. Note that all operations 
from $\Z/p$ to $\Z$ and from $\Z$ to $\Z/p$ are built from combinations of $\rho_p$ or $\beta$ 
with Steenrod powers (or squares for $p=2$).

\medskip
As differential cohomology is built out of integral cohomology and differential form data, 
 cohomology operations in both of these settings are essential for our
 construction of the refined cohomology operations. Somewhat surprisingly, we found that 
 neither has been studied to the extent that one might expect from
 such classical notions.

\paragraph{Operations from $\Z$ to $\Z$.}  
Integral cohomology operations $K(\Z, n) \to K(\Z, m)$ have been studied starting with Cartan
\cite{C}. The algebraic structure has been investigated in \cite{May} \cite{Koc} 
\cite{P1}. However, there does not seem to be a complete characterization, at least in the 
unstable case, and explicit calculations do not seem to be available in all cases. 
An exception is perhaps \cite{P2},  where the Leray-Serre 
spectral sequence for the path-loop fibration $K(\Z, n) \to PK(\Z, n+1) \to K(\Z, n+1)$ is 
used to calculate the groups $H^m(K(\Z, n); \Z)$ for $2\leq n \leq 7$ and $2 \leq m \leq 13$,
which should be useful for applications.  So, aside from arriving at $\Z$-operations 
via the Steenrod algebra, not much seems to be known. 

\paragraph{Operations from $\Omega^n$ to $\Omega^m$.}
Unlike all the above, these operations are not at the level of cohomology, but rather 
occur at the level of differential forms. 
For compact manifolds, linear operations on differential forms $\Omega^n(X) \to \Omega^m(X)$ 
which commute with diffeomorphisms have been considered 
in Palais \cite{Pa} from the point of view of functional analysis. This was extended to 
the noncompact case in \cite{J}. 
More general operations in a much broader context are studied in \cite{KMS},
but the operations relevant to us are still linear (and we are interested in nonlinear ones as well); 
there it is shown that all operations that raise the form degree by one are multiples of the exterior derivative,
and linearity follows from naturality. More recently, operations (both linear and nonlinear) acting on 1-forms (connections) 
were considered in \cite{FH}, and generalized to differential forms of all degrees in \cite{NS}.
We will make use of this for our construction of cohomology operations on
closed differential forms $\Omega^*_{\rm cl}$ in stacks.

\medskip
Cohomology operations need not be homomorphisms. Indeed, 
the power map $H^n(X; G) \to H^{2n}(X; G)$ is a cohomology 
operation by naturally, but is obviously not a homomorphism.
However, this map becomes a homomorphism when $G=\Z/p$,
and it is the example of a top degree Steenrod operations.

\medskip
One might wonder whether there is anything at all to be gained by considering
differential cohomology operations, as after all we are considering operations that 
emerge from $\Z/p$ coefficients, while the additional data in differential cohomology 
is that of de Rham forms. On the other hand, one would think that there must be some effect of 
differential refinement on the cohomology operations, as after all 
there are differential refinements of characteristic cohomology classes that 
led to a considerable amount of utility and applications (see \cite{HS} \cite{SSS3} \cite{FSSt} \cite{FSS1}
\cite{FSS2} \cite{E8}).  It turns out that the general result will fall somewhere in between. 
For instance, we find that
\begin{itemize} 
\item Even Steenrod squares cannot be differentially refined. 
\item Odd Steenrod squares refine as 
 $ \widehat{Sq}^{2m+1}=j\Gamma_2 Sq^{2m} \rho_2 I$,
where
\begin{itemize}
\item $I: \widehat{H}^*(-; \Z) \to H^*(-; \Z)$ the integration map, corresponding to `unrefinement'. 
\item $\Gamma_p$ is induced by a representation $\ZZ/p\into U(1)$ as the roots of unity
\item $j$ is the flat inclusion $H^*(-, U(1)) \hookrightarrow \widehat{H}^*(-; \Z)$, i.e. inclusion 
of flat bundles into bundles with connection.
\end{itemize} 
\end{itemize}
Note that the special case of Steenrod squares in degree one less than the top degree 
were considered in Gomi \cite{Go} (see also \cite{Bun1} Sec. 3.4) and related to the Deligne-Beilinson cup product. Thus, a portion of our work can be viewed as a generalization of this relationship to all degrees.
 
\medskip
The paper is organized as follows.
The first two subsections of Sec. \ref{Sec Form} are meant to 
give a directed overview of the two main ingredients that we aim 
to coherently merge together, namely Steenrod operations and stacks. 
In Sec. \ref{Sec Class} we recall the definition of 
Steenrod squares and Steenrod powers from two points of view:
via (co)chains and via symmetric group actions. We present these 
in such a way that helps the reader conceptually follow the constructions 
in later sections. Here we also recall the integral lifts of Steenrod squares,
which are needed for differential refinements. 
In Sec. \ref{Sec Stack} we set up the machinery of stacks, adapted to 
the context of differential cohomology,
that we 
need in order to formulate the differential refinements. 
The main general results are presented in Sec. 
\ref{Sec Diff},  where we present the characterization theorem 
(Theorem \ref{characterization theorem}) of general differential 
cohomology operations. This requires an interplay between 
integral cohomology operations and operations on differential forms. 

\medskip
We apply this formulation to the Steenrod operations 
in Sec. \ref{Sec Sq}, where the even case is given in 
Proposition  \ref{Prop even}, while the odd case is 
given in Corollary \ref{Cor Sq odd}. 
In Sec. \ref{Sec DB} we investigate whether or 
not the differential Steenrod squares are related to the homotopy commutativity 
of the Deligne-Beilinson cup product \cite{De} \cite{Be} (see also \cite{Bry}),
refining the classical point of view on the Steenrod squares presented in Sec. \ref{Sec Class}.
This leads to a generalization of \cite{Go} to Steenrod square of all degrees
and at the level of stacks.
 Along the way we prove a Kunneth decomposition for differential cohomology 
 (Proposition \ref{Prop Kunneth}) which should be interesting in its own right as a general 
  computational tool. 
The properties of the refined Steerod operations are given in 
 Sec. \ref{Sec Prop}. Most of the properties of the classical operations 
 continue to hold with the exception of the identity and the Cartan formula, both
 of which can be traced to the fact that the even Steenrod squares do not refine. 
 Finally we present applications of our operations in Sec. 
\ref{Sec App}, where we also introduce a special version of stability for the operations
(Proposition \ref{Prop stable}) and end with tantalizing applications to physics via
higher geometry.

\section{Formulation of primary operations in differential cohomology} 
\label{Sec Form} 

\subsection{Classical cohomology operations via (co)chains and via symmetric group actions} 
\label{Sec Class}

The material in this section is standard (\cite{MT} \cite{St}), but we include it as it helps 
in the conceptual understanding of our constructions later, due to the similarity of the
structure involved. 

\medskip
Steenrod squares are initially meant to square a class $x$, $|x|=n$ of the same degree as the operation,
i.e. $Sq^i(x)=x^2$ when $i=n$. The lower degree operations $Sq^i$, $0 \leq i \leq n-1$, 
measure in a precise way the extent to which homotopy commutativity of the cup product
deviates from strict commutativity. The cup product is not (graded)-commutative
at the chain level, and the obstruction is measured by the Steenrod operations.
We will closely follow the presentation in \cite{Bru} (Ch. 3) for an illuminating illustration of how the 
commutativity vs. homotopy commutativity arise.  

\medskip
The diagonal map $\Delta: X\to X\times X$ leads to a strictly commutative triangle at the chain level 
$$
\xymatrix{
C_*(X)  \ar[rr]^{\Delta_*} 
\ar[drr]_{\Delta_*} &&  
C_*(X \times X)
\\
&&
 C_*(X \times X)\;,
\ar[u]_{\tau_*}
}
$$ 
where  $\tau_*$ is the map induced from the exchange 
map $\tau: X \times X \to X \times X$
given by $\tau(x, y)=(y,x)$.
Now the Alexander-Whitney map 
$
\xymatrix{
C_*(X \times X) 
\ar[r]^-{\rm AW}_-{\cong} 
&
C_*(X) \otimes C_*(X)
}
$,
defines an equivalence of chain complexes but is only homotopy commutative and not 
strictly so (see \cite{Mac}). Dualizing to cochains, we can define the define the cup product operation by $(a \cup b)(x)=(a \otimes b)({\rm AW}(\Delta_*(x)))$. The effect of the homotopy commutativity of the Alexander-Whitney map propagates to the cup product and we get a diagram
$$
\xymatrix{
C^*(X) \otimes C^*(X) 
\ar[rr]^{\cup} 
\ar[d]^{\tau}
&&
C^*(X) \\
C^*(X) \otimes C^*(X) 
\ar[rru]^\cup
&&
}
$$
which is only homotopy commutative and not strictly commutative. 
This then induces the commutative diagram at the level of mod 2 cohomology
$$
\xymatrix{
H^*(X; \Z_2) && H^*(C^*(X) \otimes C^*(X); \Z_2) \ar[ll]_<<<<<<<{\mu}
\\
&& H^*(C^*(X) \otimes C^*(X); \Z_2)\;. 
\ar[u]_{\tau^*}
\ar[llu]^<<<<<<<<<<<<<<<<<<{\mu}
}
$$
 This is not yet a multiplication, for which one needs the 
K\"unneth isomorphim $H^*(X\times X;\ZZ/2)\simeq H^*(C^*(X) \otimes C^*(X); \Z_2) \cong H^*(X; \Z_2) \otimes H^*(X; \Z_2)$.

\medskip
Having homotopy commutativity then allows for a lot of structure, arising from 
the chain homotopies and then homotopies on these, all the way up until 
the dimensions are exhausted. At the first level, one gets a chain homotopy 
$\cup_1$ between $\cup \tau$ and $\cup$ corresponding to the homotopy 
$\cup \tau \simeq \cup$, such that 
$ b \cup a - a \cup b= \cup \tau (a \otimes b) - \cup (a \times b)= d \cup_1(a \otimes b) + \cup_1 d (a \otimes b)$.
Now considering the case $b=a$, we get 
$a \otimes a - a \otimes a=0=d \cup_1 (a \otimes a) + \cup_1 d (a \otimes a)$. 
If, furthermore, $a$ is taken to be a cocycle, i.e. $da=0$, then $d(a \otimes a)=0$ as well. 
Then we are left with one factor, $d \cup_1 (a \otimes a)=0$, which defines  a cohomology
class  at the next lower level 
$$
Sq^{n-1} (a):=a \cup_1 a \in H^{2n-1}(X; \Z_2)\;.
$$

\medskip
The lower Steenrod squares are obtained as the higher chain homotopies, obtained
by iterating the above process to $\cup_{i+1}: \cup_i \tau \simeq \cup_i$ for each $i \geq 0$. 
These give the remaining Steenrod squares
$$
Sq^{n-i}(a):= a \cup_i a \in H^{2n-i}(X; \Z_2)\;.
$$
The process stops at after $n$ steps, when we reach $Sq^0$, which is 
the identity.

\medskip
 Steenrod powers $P^i$, at a prime $p$, work similarly by replacing $\tau$ with the cyclic 
permutation operation on the product of $p$-fold copies
of the space $X$. This then gets translated analogously to a power map on 
cohomology classes.

\medskip
Note that one does not necessarily need to deal with chain complexes in order to construct 
the Steenrod operations. 
In fact, there is an analogous construction in topological spaces, i.e in the category $\mathscr{T}{\rm op}$,
 which makes use of the representability of 
cohomology via Eilenberg-MacLane spaces (see e.g. \cite{Hat}).
As in the above description, one begins with a homotopy commutative diagram
$$
\xymatrix{
X 
\ar[rr]^{\Delta} \ar[rrd]_{\Delta}&& 
X \times X \ar[rr]&& K(\ZZ/2,n)\times K(\ZZ/2,n)\ar[r]^-{\cup} & K(\ZZ/2,n)\\
&& X \times X \ar[rru]
\ar[u]^{\tau} &
}
$$ 
describing the homotopy commutativity of the cup product. Since we are concerned 
with the square cup product, we choose the maps to the Eilenberg-MacLane spaces 
to be given by the same map on each factor in the product, that is,
$$
x\times x:X\times X\to K(\ZZ/2,n)\times K(\ZZ/2,n)\;.
$$
Then a homotopy from this map to itself represents a loop: that is a map (analogous to the chain homotopy 
$\cup_1$ in the above approach)
$$
h:X\times X\times S^1\to K(\ZZ/2,n)\;,
$$
defined via the equivalence ${\rm Map}(S^1, {\rm Map} (X \times X, K(\Z/2, n))\simeq \map(S^1\times X\times X, K(\ZZ/2))$.
Choosing $h$ to be nontrivial, one can iterate this process and extend this map to the infinite-dimensional 
 sphere $S^\infty$ (a process that is analogous to choosing the higher homotopies in the above approach). 
 Using the symmetry of the cup product, one can choose this map in such a way that it commutes 
 with the $\ZZ/2$-action on $X\times X$ (given by transposing the factors) and the $\ZZ/2$-action on $S^{\infty}$ 
 (given by the antipodal action). The map $h$ then descends to a map on the quotient of the diagonal action
$$
h:X\times X\times_{\ZZ/2} S^{\infty}\to K(\ZZ/2,n)\;.
$$
Taking the trivial $\ZZ/2$-action on $X$, we see that the diagonal map $\Delta: X\to X\times X$ is equivariant with respect to 
this action. The same map then induces a map on the corresponding homotopy orbits and the entire construction can be represented diagrammatically as follows. 
The map $h$ is the universal map filling the homotopy commutative diagram
\begin{center}
\begin{tikzpicture}
\label{homotopy diagram}
\matrix (m) [matrix of math nodes, column sep=3em, row sep=3em]{
{X}\times \RR P^{\infty} & {X} \times {X} \times_{\ZZ/2} S^{\infty}&
\\
{X} & {X} \times {X} & K(\ZZ/2,2n)\;.
\\
};

\path[->] 
(m-1-1) edge node[above] {\footnotesize $Q(\Delta)$} (m-1-2)
(m-2-1) edge (m-1-1)
(m-2-2) edge (m-1-2)
(m-2-1) edge node[above] {\footnotesize $\Delta$} (m-2-2)
(m-2-2) edge node[above] {\footnotesize $\cup$} (m-2-3)
(m-1-2) edge[dashed] node[above] {\footnotesize $h$} (m-2-3);
\end{tikzpicture}
\end{center}
In fact, this diagram can be summarized in terms of $(\infty,1)$-colimits by identifying the map $h$ as arising from the universal property of $(\infty,1)$-colimits. Although the construction of the map $h$ is quite classical, its interpretation as the universal map filling the homotopy commutative diagram seems to be a new idea; one which we consider to have a distinct conceptual advantage. We will make use of this type of construction explicitly later in this paper (see Proposition \ref{Prop cup h}). 

\medskip
To get to the construction of the Steenrod squares, we notice that the composite map 
$hQ(\Delta):X\times \RR P^{\infty}\to K(\ZZ/2,2n)$ can be identified with an element in 
$Sq\in H^*(X\times \RR P^{\infty};\ZZ/2)$. The Ku\"nneth formula allows us to expand this element as a polynomial in powers of the first Stiefel-Whitney class.
Define the Steenrod squares to be the coefficients of the polynomial 
\(
\label{Sq poly}
Sq=Sq^n+Sq^{n-1}\otimes w_1+Sq^{n-2}\otimes w_1^2+\hdots+Sq^0\otimes w_1^n\;.
\)
One can think of the various powers of $w_1$ as indexing the cells of a CW structure, built from the transposition map. The coefficients attach homotopies involved in the cup product to these cells.

\medskip
The Steenrod squares satisfy several desirable properties \cite{St} \cite{MT} \cite{MiS} \cite{Hat}. 
To avoid unnecessary 
redundancies, we will not record these. The differential refinement of these
cohomology operations will satisfy some of the same properties, generalized
appropriately, but with marked differences. Hence, we prefer that these classical  
properties be deduced from the refined one (see Sec. \ref{Sec Prop}).

\medskip
We need to discuss integral cohomology operations on the
path to arriving at differential cohomology operations. 
In particular, we would like to consider integral lifts of the Steenrod squares, that is we seek a diagram
$$
\xymatrix{
H^\bullet(-; \Z) \ar[rr]^{Sq^j_\Z} \ar[d]^{\rho_2} && H^{\bullet+j}(-; \Z) \ar[d]^{\rho_2}
\\
H^\ast(-; \Z/2) \ar[rr]^{Sq^i} && H^{\ast+i}(-; \Z/2)
}
$$
where $\bullet$, $\ast$, $i$, and $j$ are degrees to be determined. 
In the case where $i$ is even, it is impossible to find such a lift. This can be immediately deduced using the long exact Bockstein sequence corresponding to the short exact sequence
$
\ZZ\overset{\times 2}{\longrightarrow} \ZZ \overset{\rho_2}{\longrightarrow} \ZZ/2
$.
Indeed, if such an integral refinement existed, we would have an integral operation $\theta$ such that $\rho_2(\theta)=Sq^{2i}$. But by the long exact sequence, this would imply $\beta(Sq^{2i})=0$
 (where $\beta$ is the Bockstein corresponding to the above sequence). However, this is not the case as
$
 \rho_2\beta(Sq^{2i})=Sq^{2i+1}
 $.
In contrast, for the odd Steenrod squares, the above formula allows us to define an integral refinement
\(
Sq_{\ZZ}^{2i+1}:=\beta(Sq^{2i})\;.
 \label{Integral Sq} 
\)
Our results on differential Steenrod squares will end up having a similar pattern. 
The even case will be presented in Proposition  \ref{Prop even}, while the odd case is 
given in Corollary \ref{Cor Sq odd}.

\subsection{Stacks associated to differential cohomology}
\label{Sec Stack}

In this section, we review some some of the stacks with which we will be working, and highlight some of their 
properties that are useful for us. 
A more complete study can be found in \cite{FSSt}\cite{SSS3} \cite{FSS1}\cite{FSS2}.

\medskip
To start, in order to properly account for our stacks we have the following. 
\begin{definition}
Let $\cartsp$ denote the category whose objects are open subsets $U\subset \RR^n$ that are diffeomorphic to $\RR^n$, for some $n\in \NN$, and whose morphisms are smooth maps $
f:U\to V
$
. This category admits the structure of a Grothendieck site, with topology generated by good open covers (contractible finite intersections). 

We denote $\infty$-category of $\infty$-presheaves on $\cartsp$ by ${\rm PSh}_{\infty}(\cartsp)$. We denote the $\infty$-category of $\infty$-sheaves on $\cartsp$ by ${\rm Sh}_\infty(\cartsp)$. 
\end{definition}

\begin{remark}
The category of $\infty$-sheaves on $\cartsp$ admits a presentation by model categories. More precisely, taking the homotopy coherent nerve on fibrant/cofibrant objects in the {\v C}ech local, projective (or injective) simplicial model structure on simplicial presheaves on $\cartsp$ yields a presentation. This allows us to present $\infty$-sheaves via projectively fibrant simplicial presheaves that satisfy homotopy descent.
\end{remark}

\begin{definition}
A local object $\mathcal{G}\in {\rm Sh}_\infty(\cartsp)$ is called a \emph{smooth stack} (or higher stack). In the presentation by simplicial presheaves, such an object can be identified with a functor
$$
{\cal G}:\cartsp\to \mathscr{K}{\rm an}\subset \sset\;,
$$
that satisfies homotopy descent with respect to {\v C}ech covers (see \cite{Dug}\cite{Lur}\cite{Urs} for review)\;.
\end{definition}

An important functor that allows for passage between simplicial sets and chain complexes 
in positive degrees is the Dold-Kan functor 
\(
DK: \chp \longrightarrow \mathscr{K}{\rm an}\subset\sset\;.
\label{eq DK}
\)
The stack that we will use most frequently in this paper is the moduli stack of $n$-bundles (or gerbes) 
with connection: $\BB^nU(1)_{\nabla}$. This stack arises as a pullback of stacks \footnote{Whenever we say pullback, pushout, limit or colimit, we mean these operations in the $(\infty,1)$-sense. With an appropriate choice of model structure, these can be thought of as the \emph{homotopy} pullback.}
\(\label{homotopy pullback1}
\xymatrix{
\BB^nU(1)_{\nabla}\ar[r]\ar[d] & {\Omega}^{n+1}_{\rm cl}\ar[d]
\\
\BB^{n+1}\ZZ\ar[r]^-{i} & \flat_{\rm dR}\BB^{n+1}U(1)
}
\)
If we forget about the connection on the these $n$-bundles, we obtain the bare moduli stack of $n$-gerbes $\BB^{n}U(1)$. Explicitly, this stack is obtained by applying the Dold-Kan functor 
to the sheaf of chain complexes 
$C^{\infty}(-,U(1))[n]$: the sheaf of smooth $U(1)$-valued functions in degree $n$. The other stacks related to $\BB^{n}U(1)_{\nabla}$ in the above diagram are defined via the Dold-Kan 
correspondence \eqref{eq DK} as follows (see \cite{FSSt} \cite{FSS1} \cite{FSS2} \cite{Urs}):
\begin{itemize}
\item The stack $\BB^{n+1}\ZZ$ is defined to be the smooth stack obtained by applying the Dold-Kan functor to the sheaf of chain complexes $\underline{\ZZ}[n]$, with the sheaf of locally constant $\ZZ$-valued functions in degree $n$. 

\item The stack representing the truncated de Rham complex $\flat_{\rm dR}\BB^{n}U(1)$ is obtained by applying Dold-Kan to the truncated de Rham sheaf of chain complexes 
$$\Omega^{\leq n}_{\rm cl}:=[\hdots 0 \hdots \Omega^{0}\to \Omega^{1} \to \hdots \Omega^{n}_{\rm cl}]\;.$$

\item The stack of closed $n$-forms ${\Omega}^{n}_{\rm cl}$ is defined 
to be the stack obtained by applying Dold-Kan to the sheaf of closed $n$-forms.
This is also discussed in \cite{HS} and \cite{Bun1} (Problem 4.42).  

\end{itemize} 

The differential cohomology diagram \cite{SSu}
 lifts to a diagram of stacks  \cite{Bun1} \cite{Urs}
\(
\label{stackdiamond}
\begin{tikzpicture}
\matrix (m) [matrix of math nodes, column sep=3em, row sep=3em, nodes={anchor=center}]{
& {\Omega}^{\leq n} & & {\Omega}^{n+1}_{\rm cl}&
\\
 \flat_{\rm dR}\BB^{n-1}U(1) & &  \BB^{n}U(1)_{\nabla} & & \flat_{\rm dR} \BB^{n+1}U(1) 
\\
& \flat\BB^{n}U(1) & & \BB^{n+1}\ZZ &
\\
};

\path[-stealth]
(m-1-2) edge node[above] {\small ${d}$} (m-1-4)
(m-2-1) edge (m-1-2)
(m-1-2) edge node[above right] {\small ${a}$} (m-2-3)
(m-2-3) edge[->>] node[above right] {\small ${I}$} (m-3-4)
(m-2-1) edge (m-3-2)
(m-3-2) edge node[above right] {\small ${\beta}$}(m-3-4) 
(m-3-4) edge (m-2-5)
(m-3-2) edge (m-2-3)
(m-2-3) edge[->>] node[below right] {\small ${R}$} (m-1-4)
(m-1-4) edge (m-2-5);     
;
\end{tikzpicture}
\)
where the diagonals are fiber sequences. Moreover, the maps $a$, $I$ and $R$ induce homomorphisms in cohomology. In fact, the stacks surrounding the center are related to the center in a functorial way. The functors which produce these surrounding stacks are part of an $\infty$-adjunction called a \emph{cohesive adjunction} \cite{BNV} \cite{Urs}. 

\medskip
It is shown in \cite{Urs} that 
the $\infty$-category of smooth higher stacks ${\rm Sh}_{\infty}(\cartsp)$ admits a quadruple $\infty$-categorical adjunction 
$(\Pi \dashv {\rm disc} \dashv \Gamma \dashv {\rm codisc})$ 
\(
\xymatrix{
{\rm Sh}_{\infty}(\cartsp)\ar@<-.5em>[rr]^<<<<<<<<<<{\Gamma} \ar@<1.5em>[rr]^<<<<<<<<<<{\Pi} && \sset \ar@<1.5em>[ll]_<<<<<<<<<<{\rm codisc} \ar@<-.5em>[ll]_<<<<<<<<<<{\rm disc} 
}
\label{cohesion}
\)
where the fundamental groupoid functor $\Pi$ preserves finite $\infty$-limits,
 and the discretization and co-discretization functors ${\rm disc}$ and ${\rm codisc}$ are fully faithful. In our context the quadruple adjunction is presented by Quillen functors, where on the left we take the {\v C}ech local projective model structure.  
 
 One implication of this is that $\sset$ embeds into ${\rm Sh}_{\infty}(\cartsp)$ in two different ways as a 
reflexive $\infty$-subcategory. From the reflectors we can produce two monads and one comonad defined as follows:
\(
{\bf \Pi}:=\Pi\circ {\rm disc},\qquad \quad
 \flat:={\rm disc}\circ \Gamma, \qquad \quad
  \sharp:={\rm codisc}\circ \Gamma\;.
 \label{monad-adj}
 \)
These monads fit into a triple adjunction $({\bf \Pi}\dashv \flat \dashv \sharp)$ which is called a \emph{cohesive} adjunction. We remark that the adjoint functors can be presented by derived functors in the appropriate Quillen model structure (either local injective or local projective) on simplicial presheaves (see \cite{Urs} for details). 

\begin{remark}
Each monad in the cohesive adjunction picks out a different part of the nature of a smooth stack. This nature is perhaps best exemplified by how the adjoints behave on smooth manifolds (viewed as stacks). More precisely, if $M$ is a smooth manifold
then, for instance, 
\begin{enumerate}[label=(\roman*)]
\item 
the comonad $\flat$ (flat) takes the underlying set of points of the manifold and then embeds this set back into stacks as a discrete object. This functor therefore misses the smooth structure of the manifold and treats it instead as a discrete object. 

\item The monad ${\bf \Pi}$ essentially takes the singular nerve of the manifold using \emph{smooth} paths and higher smooth simplices on the manifold. It therefore retains the geometry of the manifold and ``knows" that the points of the manifold ought to be connected together in a smooth way. 
\end{enumerate}
\end{remark} 
As discussed in \cite{Urs}, the bare stack $\BB^n\ZZ$ representing integral cohomology is equivalent to ${\bf \Pi}\BB^nU(1)_{\nabla}$, while the discrete stack $\BB U(1)^{\delta}$ is equivalent to the moduli stack of flat bundles $\flat \BB^nU(1)_{\nabla}$. Thus, one can rewrite the diamond diagram \eqref{stackdiamond} using only the monads ${\bf \Pi}$ and $\flat$. In this context the ``unrefinement map" $I$ arises as the unit of the monad, i.e. a natural transformation 
\(
I:{\rm id}\to {\bf \Pi}\;,
\label{eq I}
\)
where ${\rm id}$ is the identity functor. 
In the stable case, the characterization of differential cohomology theories using these monads is due to \cite{BNV}.

\medskip
We now prove a few properties which we will use later.
\begin{proposition}
\label{Prop-disc}
Let $A$ be a discrete abelian stack. That is $A\simeq {\rm disc}(B)$ form some $B\in {\rm s}\mathscr{A}{\rm b}$. Then
$$\flat \BB^nA\simeq \BB^n A\;.$$
\end{proposition}
\theproof
Since $\Gamma$ and ${\rm disc}$ are right adjoints, they commute with looping, as this operation 
is an example of an $\infty$-limit. Consequently they also commute with delooping, as this is defined
using looping. It follows by definition (see \eqref{monad-adj}) that $\flat$ commutes with delooping. 
Then we have
\bea 
\flat \BB^nA &\simeq& \BB^n \flat A 
\\
&=& \BB^n ({\rm disc}\circ \Gamma)(A)
\\
&\simeq& \BB^n ({\rm disc}\circ \Gamma \circ {\rm disc}(B))
 \\
 &\simeq & \BB^n({\rm disc}(B))
 \\
 &\simeq & \BB^n A\;.
 \eea
 In the step before last we used the fact that $\Gamma \circ {\rm disc}={\rm id}$, because 
 ${\rm disc}$ is given by the stackification of the constant functor $U \mapsto A$, and then 
 $\Gamma$ evaluates that at a point $U=\R^0$, giving back the original object $A$. 
 \endofproof

The moduli stacks $\BB^nU(1)_{\nabla}$, as $n$ varies, represent differential cohomology in the sense that the functor
$$\pi_0\map(-,\BB^nU(1)_{\nabla}):{\rm Sh}_{\infty}(\cartsp)\to \ab\;$$
assigns to every smooth manifold $X$ (embedded in stacks via the sheaf of smooth plots $C^{\infty}(-,X)$) the differential cohomology group $\widehat{H}^{n+1}(X;\ZZ)$. This follows almost immediately from the presentation of ordinary differential cohomology as Deligne cohomology, along with the Dold-Kan correspondence \cite{FSSt} \cite{Urs}. 
The unrefinement morphism $I:{\rm id}\to {\bf \Pi}$ then induces a natural transformation 
$$
 I_*: \pi_0\map(-,\BB^nU(1)_{\nabla}) \to \pi_0\map(-,{\bf \Pi}\BB^nU(1)_{\nabla}) 
\simeq  \pi_0\map(-,\BB^{n+1}\ZZ)\;.
$$
We can equivalently write this map as a map
$$
I_*:\widehat{H}^{n+1}(-;\ZZ) 
\to H^{n+1}(-;\ZZ)\;,
$$
which is the familiar ``integration map" which relates differential cohomology 
to its underlying cohomology theory.

\medskip
We will frequently need to use representability in our calculations and indeed
we are able to pass from the stacks to underlying theories somewhat seamlessly. 
For this reason, we remind the reader of the various theories that are represented 
by the stacks in the refined diamond diagram \eqref{stackdiamond} 
(see \cite{FSS1} \cite{FSS2} \cite{E8} \cite{Urs} for details). Below we use the notation $\map(-,-)$ for the derived simplicial hom.

\begin{itemize}
\item The stack $\BB^{n+1}\ZZ$ represents ordinary integral cohomology in degree $n+1$. 
That is, we have a natural isomorphism 
$$
\pi_0\map(-,\BB^{n+1} \ZZ)\simeq H^{n+1}(-;\ZZ)\;.
$$
\item The stack $\flat \BB^n U(1)$ represents cohomology with $U(1)$-coefficients:
$$
\pi_0\map(-,\flat \BB^n U(1))\simeq H^{n}(-;U(1))\;.
$$
\item The stack $\flat_{\rm dR}\BB^{n}U(1)$ represents de Rham cohomology in degree $n+1$,
$$
\pi_0\map(-,\flat_{\rm dR}\BB^nU(1))\simeq H_{\rm dR}^{n+1}(-)\;.
$$
Equivalently, by de Rham's theorem, this stack also represents cohomology with real coefficients $H^{n+1}(-;\RR)$.
\item The stack $\Omega^n_{\rm cl}$ represents the sheaf of closed differential $n$-forms:
$$\pi_0\map(-,\Omega^{n}_{\rm cl})\simeq \Omega^{n}_{\rm cl}(-)\;.$$
\end{itemize}
In some of the proofs, we will use the properties of discrete stacks when calculating homotopy classes of maps. More precisely, we have an adjunction 
\(\label{nerve}
\xymatrix{
\sset \ar@<.2cm>[rr]^{\pi_0} &&  \set \ar@<.2cm>[ll]_{\rm sk_0}
}\;,
\)
with $\pi_0$ takes the connected components of the simplicial set. The right adjoint ${\rm sk}_0$ simply embedds a set as a discrete stack. To illustrate how one uses the adjunction in practice, observe that both a manifold $X$
and the sheaf of closed $n$-forms $\Omega^n_{\rm cl}$ are represented by discrete objects in the category of stacks. 
Then we have
$$
\pi_0\map(X, \Omega^{n}_{\rm cl})\simeq \pi_0\map(X, {\rm sk}_0(\Omega^{n}_{\rm cl}))
\simeq 
\hom(\pi_0(X),\Omega^n_{\rm cl})\simeq \Omega^n_{\rm cl}(X)\;.
$$
Here we have used the adjunction between ${\rm sk}_0$ and $\pi_0$, passing to the category of sheaves, and 
finally used the Yoneda lemma.

\subsection{General differential Cohomology operations} 
\label{Sec Diff} 

We consider differential cohomology operations from a general point of view,
and then specialize in later sections. We will consider these operations in the 
context of the stacks approach to differential cohomology.

\medskip
As indicated in the Introduction, 
we will need to study morphisms of stacks
\(\label{381}
\widehat{\theta}:\mathbf{B}^nU(1)_{\nabla}\to \mathbf{B}^mU(1)_{\nabla}\;.
\) 
The goal of this section will be to establish the general properties of these maps and to provide a characterization theorem describing the general form of all such operations. In order to prove the theorem, we will need to understand how differential cohomology operations refine classical ones. The cohesive adjoints \eqref{cohesion} will be extremely useful in our discussion on this point. Essentially, this boils down to the fact that these functors pick out different aspects of the moduli stack $\BB^nU(1)_{\nabla}$. Hence, when studying the maps \eqref{381}, we can use these functors to isolate various parts of the source or target stack (depending on the situation). We can then use the properties of these functors to arrive at isomorphisms between hom sets in the homotopy category which would otherwise require a lot of work to establish.

\medskip
In what follows, we will continue to denote the differential cohomology group of a differentiable manifold $X$ in degree $n$ by $\widehat{H}^n(X;\ZZ)$ and identify it with the contravariant functor
$$
\widehat{H}^n(-;\ZZ):=\pi_0\map(-,\BB^nU(1)_{\nabla})\;,
$$
restricted to the subcategory of smooth manifolds.

\begin{definition}
A \emph{differential cohomology operation} is a natural transformation of functors
$$
\widehat{\theta}:\widehat{H}^n(-;\ZZ) \to \widehat{H}^m(-;\ZZ)
$$
or, equivalently, a homotopy class of maps between stacks
$$
\widehat{\theta}:\BB^nU(1)_{\nabla}\to \BB^mU(1)_{\nabla}\;.
$$
\label{Def coh op}
\end{definition}
At this stage, we can already see one of the advantages provided by the stacky approach to differential cohomology operations. 
This is, we can describe these operations as elements in the set
$\pi_0\map(\BB^nU(1)_{\nabla},\BB^mU(1)_{\nabla})$, just as integral cohomology operations are elements in $H^m(K(\ZZ,n);\ZZ)\simeq \pi_0\map(K(\ZZ,n),K(\ZZ,m))$. This allows us to do constructions at the universal level, 
which is not possible without the use of stacks. 

\medskip
Now since the stack $\BB^nU(1)_{\nabla}$ arises as the pullback \eqref{homotopy pullback1}, 
the universal property of pullbacks ensures that a map $\widehat{\theta}$ of the type \eqref{381}
is induced by a homotopy commutative diagram involving operations $\tau$ on closed
differential forms, $\alpha$ on de Rham cohomology, and $\theta$ on singular cohomology
\(\label{homotopy pullback}
\xymatrix{
&
{\Omega}^{n+1}_{\rm cl}
 \ar[rr]^{\tau}
 \ar[d]^c
 &&
{\Omega}^{m+1}_{\rm cl}
\ar[d]^c
\\
\BB^n U(1)_{\nabla} \ar[ur]^R \ar[dr]_I
&
\flat_{\rm dR}\BB^{n}U(1)
\ar[rr]^{\alpha}
&&
\flat_{\rm dR}\BB^{m}U(1)
\\
& \BB^{n+1}\ZZ
\ar[rr]^{\theta}
\ar[u]_{i}
&&
\BB^{m+1}\ZZ\;.
\ar[u]_{i}
}
\)
Hence, we see that a triple $(\theta,\alpha,\tau)$ which makes the above diagram commute up to a choice of 2-morphism immediately induces a differential cohomology operation $\widehat{\theta}$. This point of view emphasizes that differential cohomology operations are really a compatible combination of operations on differential forms and operations on integral cohomology. Furthermore, these operations are required to be homotopic in the de Rham stack $\flat_{\rm dR}\BB^{m}U(1)$. 

\medskip 
If a differential cohomology operation $\widehat{\theta}$ is induced from such a triple, we
 say that $\widehat{\theta}$ \emph{refines} the integral cohomology operation 
$\theta$ and the operation $\tau$, on differential forms. On general abstract grounds, the 
converse of the above statement 
may not be true. That is, every morphism $\widehat{\theta}$ as in Definition \ref{Def coh op} 
 need not be induced by such a triple. 
However, because of the special nature of the particular pullback involving 
$\BB^nU(1)_{\nabla}$, we now see that this is indeed the case.

\begin{proposition}
\label{refinement}
Every differential cohomology operation $\widehat{\theta}$ refines an integral operation 
$$
\theta:\BB^{n+1}\ZZ\to \BB^{m+1}\ZZ
$$
and an operation on 
forms
$$
\tau:\Omega^{n+1}_{\rm cl}\to \Omega^{m+1}_{\rm cl}
$$
which fit into a homotopy commuting triple $(\theta,\alpha,\tau)$ as in \eqref{homotopy pullback}. Moreover, at the level of homotopy classes, we have
$$
\tau R=R\widehat{\theta}\ \ \text{and} \ \ \  \theta I=I\widehat{\theta}\;.
$$
Diagrammatically, we have homotopy commutativity 
\(
\label{house}
\xymatrixcolsep{4pc}\xymatrixrowsep{3.9pc}
\xymatrix{
& \BB^n U(1)_{\nabla} 
\ar[r]^{\widehat{\theta}}
\ar@/_1.2pc/[dl]_{R}
\ar@{.>}@/^.8pc/[dddl]_{I}
& 
\BB^m U(1)_{\nabla} 
\ar@/_1.2pc/[dl]_{R}
\ar@/^1.5pc/[dddl]^{I}
\\
{\Omega}^{n+1}_{\rm cl}
\ar[r]^\tau 
\ar[d]^c
&
{\Omega}^{m+1}_{\rm cl}
\ar[d]^c
\\
 \flat_{\rm dR}\BB^{n+1}U(1)
\ar[r]_{~~\alpha}
&
\flat_{\rm dR}\BB^{m+1}U(1)
\\
\BB^{n+1}\ZZ
\ar[r]^{\theta}
\ar[u]_-{i}
&
\BB^{m+1}\ZZ\;.
\ar[u]_{i}
}
\)
\end{proposition}
\theproof
Let $\widehat{\theta}$ be a differential cohomology operation. Then $I\widehat{\theta}:\BB^nU(1)_{\nabla}\to \BB^{m+1}\ZZ$. Since $\ZZ$ is discrete, we have $\BB^{n+1}\ZZ\simeq \flat \BB^{n+1}\ZZ$. This fact, along with the cohesive adjunction, imply that we have an equivalence
\bea
\map(\BB^nU(1)_{\nabla},\BB^{m+1}\ZZ) &\simeq& \map(\BB^nU(1)_{\nabla},\flat\BB^{m+1}\ZZ)
\qquad \quad \quad ~({\rm by~above})
\\
&\simeq & \map({\bf \Pi}(\BB^nU(1)_{\nabla}),\BB^{m+1}\ZZ)
\qquad \quad ({\rm by~ } \eqref{monad-adj})
\\
&\simeq& \map(\BB^{n+1}\ZZ,\BB^{m+1}\ZZ)
\eea
Since the composite equivalence between the first and third line is induced by precomposition with $I$
(see the discussion around eq. \eqref{eq I}), 
we have an operation $\theta : \BB^{n+1}\ZZ\to\BB^{m+1}\ZZ$ such that $\theta I =I\widehat{\theta}$ at the level of homotopy classes.

\medskip
To prove that $\widehat{\theta}$ refines an operation on forms, we observe that $R\widehat{\theta}:\BB^nU(1)_{\nabla}\to \Omega^{m+1}_{\rm cl}\;$ and since $\Omega^{m+1}_{\rm cl}$ is a discrete object, using the adjunction \eqref{nerve}, 
we have an isomorphism
\(
\map(\BB^nU(1)_{\nabla},\Omega^{m+1}_{\rm cl})
\simeq \hom\big(\pi_0(\BB^nU(1)_{\nabla}),\Omega^{m+1}_{\rm cl}\big)
\simeq \hom(\Omega^{n}/{\rm im}(d),\Omega^{m+1}_{\rm cl})\;.
\label{eq tau0}
\)
Here we have used the isomorphisms 
\bea
\pi_0(\BB^nU(1)_{\nabla}) &\simeq&
 \pi_0(DK(\Z_{\cal D}^\infty(n+1)) 
 \\
 &\simeq&
  H_0(\Z_{\cal D}^\infty(n+1))
  \\
&\simeq& H_0
\big[
\Z \to 
 \Omega^0 \overset{d}{\longrightarrow} 
 \cdots \overset{d}{\longrightarrow} 
 \Omega^{n-1} \overset{d}{\longrightarrow} 
 \Omega^n  
\big]
\\
&\simeq&
 \Omega^n/{\rm im}(d)\;,
\eea
The exterior derivative induces an isomorphism on sheafification
(by Poincar\'e lemma)
\(
d:L(\Omega^{n}/{\rm im}(d))\to \Omega^{n+1}_{\rm cl}\;,
\label{eq dd}
\)
where $L$ is the sheafification functor. Therefore, the right hand side of eq. \eqref{eq tau0} is isomorphic to $\hom(\Omega^{n+1}_{\rm cl},\Omega^{m+1}_{\rm cl})$. The isomorphism is exactly precomposition with the curvature (eq. \eqref{eq dd}),
 and therefore there is an operation $\tau$ on forms such that $\tau R=R\widehat{\theta}$ at the level of homotopy classes. The homotopy commutativity follows from the homotopy commutativity of the pullback diagram \eqref{homotopy pullback}. 
The homotopy commutativity of diagram \eqref{house} gives a homotopy 
$c R \widehat{\theta} \to cI \widehat{\theta}$, which we can identify with $\alpha$. 
Making this identification explicit is not particularly illuminating and we leave such details to the 
interested reader. 
\endofproof

\begin{remark}
Thinking about elements of $\widehat{H}^n(X;\ZZ)$ as higher line bundles with connection, the previous proposition makes it explicit how a differential cohomology operation can be interpreted as an operation on bundles. Moreover, the curvature and underlying integral class of the resulting bundle is obtained via some singular cohomology operation and an operation on forms. 
\end{remark}

\medskip
At this stage, the reader might wish to see some examples of differential cohomology operations. Indeed, we have the following two examples which, as we will see explicitly in Lemma \ref{Lemma dR}, are essentially
the only examples which give classes with nonzero curvature. 

\begin{example} [Dixmier-Douady class]
The homotopy class of the identity morphism
$$
{\rm DD}:={\rm id}:\BB^{n}U(1)_{\nabla}\to \BB^{n}U(1)_{\nabla}
$$
is a differential cohomology operation called the (higher) Dixmier-Douady class, as this corresponds 
topologically to the fundamental cohomology class $\iota_n$ in $H^n(K(G, n); G)$. These higher classes
are amplified in \cite{FSS1} \cite{FSS2}. It is easy to see that this class refines the operations
\bea
{\mu_{n+1}}:={\rm id}&:&\Omega^{n+1}_{\rm cl}\to \Omega^{n+1}_{\rm cl}\;,
\nonumber\\
{\iota_{n+1}}:={\rm id}&:&\BB^{n+1}\ZZ\to \BB^{n+1}\ZZ\;,
\eea
where $\mu_{n+1}$ and $\iota_{n+1}$ are both the identity map, thought of as fundamental classes for the 
corresponding stacks.
The latter is perhaps familiar from the representability of singular cohomology
via Eilenberg-MacLane spaces. For the former, note that in stacks we can 
think of the stack of closed $n$-forms as representing a cohomology theory 
as well, but now in a more general sense. 
\label{DD example}
\end{example}

\begin{example} [Power operation]
Let $m$ be a positive integer. 
We will consider the Deligne-Beilinson cup product $\cup_{\rm DB}$ on stacks 
(see \cite{FSS1} \cite{FSS2} \cite{GS1}). 
Then the $m$-fold power gives a morphism of stacks
$$
{\rm DD}^m:=\underbrace{\cup_{\rm DB} \cdots \cup_{\rm DB}}_{m-{\rm times}} ~:~
\BB^{n}U(1)_{\nabla}\to \BB^{m(n+1)-1}U(1)_{\nabla}
$$
as described in \cite{FSS1}. The homotopy class of this map is, by definition, a differential cohomology 
operation. The cup product morphism refines the singular cup product and the 
wedge product of forms \cite{FSS1}\cite{GS1}. As a consequence,
 we immediately see that this operation refines
the wedge product power and the cup product power, respectively, viewed as powers of the fundamental 
classes we encountered in example \ref{DD example}. Explicitly,
\bea
\mu_{n+1}^m=\underbrace{\wedge \cdots \wedge}_{m-{\rm times}} &:&
\Omega^{n+1}_{\rm cl}\to \Omega^{m(n+1)}_{\rm cl}\;,
\nonumber\\
\iota_{n+1}^m=\underbrace{\cup \cdots \cup}_{m-{\rm times}} &:&\BB^{n+1}U(1)\to \BB^{m(n+1)}U(1)\;.
\eea
\end{example}

The remainder of this section will be devoted to proving the following main classification theorem for differential 
cohomology operations.
\begin{theorem}[Characterization theorem]
\label{characterization theorem}
Let $\widehat{\theta}$ be a differential cohomology operation. 
Then exactly one of the following holds:
\begin{enumerate}
\item $\widehat{\theta}=n{\rm DD}$, for some $n\in \ZZ$.
\item $\widehat{\theta}=n{\rm DD}^m$, for some $n\in \ZZ$.
\item $\widehat{\theta}$ factorizes as
$$
\widehat{\theta}=j\phi I\;,
$$
\end{enumerate}
where $j:\flat\BB^{m} U(1)\into \BB^mU(1)_{\nabla}$ is the flat inclusion, 
$\phi:\BB^{n+1}\ZZ\to \flat \BB^{m}U(1)$ is an operation from singular cohomology to 
cohomology with $U(1)$-coefficients, and $I$ is the canonical morphism 
$I:\BB^nU(1)_{\nabla}\to \BB^{n+1}\ZZ$.
\end{theorem}

\begin{remark} 
{\bf (i)} This theorem is analogous to cohomology of integral Eilenberg-MacLane spaces being finite or 
not, depending on the degree. 

\vspace{1mm}
\noindent {\bf (ii)} Note that the morphisms $j$, $\phi$, and $I$ can be described more classically 
along the lines of the presentation in the Introduction, and in fact generalizing those. 
Indeed, $j$ is an operation from 
$U(1)\simeq \R/\Z$-coefficients to differential cohomology (which can be viewed in a precise sense `as' the Deligne complex $\Z_{\cal D}^\infty$), 
$\phi$ is a map from $\Z$-coefficients to $U(1)$-coefficients,
and $I$ is a map from differential cohomology (as $\Z_{\D}^\infty$) to $\Z$-coefficients.  
\end{remark}

To prove the theorem, we will need to understand the rational and integral operations along with the 
operations on forms. 
 We begin with a brief recollection 
 of integral and rational operations.

\medskip
Recall that the only operations that arise rationally are the identity and the power operations. 
Indeed, the rational cohomology ring of a rational Eilenberg-MacLane space is generated by a single 
generator (see e.g. \cite{GM} Lemma 8.5) and is a $\QQ$-polynomial algebra 
or a $\QQ$-exterior algebra, depending on parity, 
\(
H^*(K(\QQ, n); \QQ)=
\left\{
\begin{array}{ll}
\QQ[\iota_{2m}], & n=2m ~{\rm even} \\
\Lambda_\QQ[\iota_{2m+1}] ,& n=2m+1~ {\rm odd},
\end{array}
\right. 
\label{eq KQn}
\)
where $\iota_q$ is the $q$th fundamental class.

\medskip
The case of singular cohomology is of course more complicated. However, the situation is made much more 
tractable by the above rational considerations. In fact, the above implies that $H^{n+q}(K(\ZZ,n);\ZZ)$ must be finite 
when $n\ndivides q$. Otherwise, the rationalization would be nonzero in these degrees, which is not the case. 
\begin{remark}
Summarizing this, along with other  properties of these groups, we have (see e.g.  \cite{C} 
\cite{Po} \cite{FFG})
\begin{enumerate}
\item $H^{n+q}(K(\ZZ,n);\ZZ)$ is finite and independent of $n$ for $0<q<n$.
\item When $n \divides q$ then this is infinite cyclic generated by powers of the fundamental class.
\item The $p$-primary part of $H^{n+q}(K(\ZZ,n);\ZZ)$ is zero for $0<q<2p-1$
\item If $n<2p-1$, the $p$-primary part of $H^{n+2p-1}(K(\ZZ,n);\ZZ)$ is cyclic of order $p$, generated 
by the operation $(\beta_pP^1_p)(u)$, where $u$ is the fundamental class, $P^1_p$ is the 1st $P$ 
operation at the prime $p$ and $\beta_p$ is the Bockstein homomorphism for the mod-$p$ sequence.
\end{enumerate}
\label{rem KZn}
\end{remark}

We now turn to the possible operations on forms, which turn out to be
 in harmony with the operations in rational cohomology.    
\begin{lemma} 
\label{Lemma dR}
Let $\mu_n:\Omega^n_{\rm cl}\to \Omega^n_{\rm cl}$ denote the identity morphism on forms (thought of as 
a fundamental class for the sheaf of closed $n$-forms). The set $\pi_0\map(\Omega^n_{\rm cl},\Omega^*_{\rm cl})$ 
forms a graded algebra and we have
$$\pi_0\map(\Omega^n_{\rm cl},\Omega^*_{\rm cl})=
\left\{
\begin{array}{ll}
\RR[\mu_{2m}], & n=2m ~{\rm even} \\
\\
\Lambda_\RR[\mu_{2m+1}] ,& n=2m+1~ {\rm odd},
\end{array}
\right. 
$$
where the asterisk $*$ on the left hand side is a grading that is determined by the powers
on the right hand side. 
\end{lemma}
\theproof
Let $f:\Omega^n\to \Omega^*$ be a natural transformation of sheaves. In \cite{NS} it was shown that any assignment of differential forms $\omega\mapsto f(\omega)$, which is natural with respect to pullback, is given by a polynomial in $\omega$ and its derivative $d\omega$. Hence, if we restrict $f$ to the sheaf of closed forms, we see that $f$ must assigns each section $\omega\in \Omega_{\rm cl}^n$, a polynomial in $\omega$. The claim is simply a restatement of this fact.
\endofproof


\medskip
The next lemma will be needed in the proof of the main theorem, Theorem \ref{characterization theorem}. 
Essentially, the lemma shows that the only differential cohomology operations detected by de Rham 
cohomology are the rational ones. The method of proof again appeals to the cohesive adjunction, which
 extracts the relevant information from the full moduli stack of $n$-bundles.
\begin{lemma} \label{de rham ops} We have
$$
\pi_0\map(\BB^nU(1)_{\nabla},\flat_{\rm dR}\BB^mU(1))=\left\{\begin{array}{ccl}
\R && {\rm if}~n=2k \divides m~{\rm or}~n=m=2k+1,
\\
\\
0 && {\rm otherwise}.
\end{array}
\right.
$$
\end{lemma}
\theproof
By de Rham's theorem, we have $\flat_{\rm dR}\BB^mU(1)\simeq \BB^{m+1}\RR$. Since $\RR$ is discrete (as a stack), we have $\flat \BB^{m+1} \RR\simeq \BB^{m+1}\RR$. By cohesion, we have equivalences
\bea
\map(\BB^nU(1)_{\nabla},\flat_{\rm dR}\BB^mU(1)) &\simeq& \map(\BB^nU(1)_{\nabla},\flat \BB^{m+1}\RR)
\qquad \qquad ({\rm by~ above~})
\\
&\simeq & \map\left({\bf \Pi}\left(\BB^nU(1)_{\nabla}\right),\BB^{m+1}\RR\right)
\qquad ({\rm by~ } \eqref{monad-adj})
\\
&\simeq & \map(\BB^{n+1}\ZZ,\BB^{m+1}\RR)
\eea
The claim then follows from the properties of Eilenberg-MacLane spaces, specifically  
the structure in \eqref{eq KQn} and the first part of Remark \ref{rem KZn}.
\endofproof

\medskip
We are now ready to prove Theorem \ref{characterization theorem}. We will find that the first case is straightforward due to the appeal to de Rham theory, while for the second case things 
become very subtle due to torsion. 

\medskip
\theproof
Let $\widehat{\theta}$ be a differential cohomology operation. By Lemma \ref{Lemma dR}, 
we have two possibilities for the corresponding operation $\tau$ on forms.

\vspace{0.1in}
 
\noindent (i) \fbox{$\tau=\lambda\mu_n^q,\ \ q\geq 0,\ \lambda\in \RR$}

\vspace{0.1in} 

First consider the case $\lambda=1$. Then $\tau$ admits at least one refinement, since ${\rm DD}^q$ refines this operation. To see that this is the only possibility, let $\hat{\theta}$ be another operation refining $\tau$. Then since $n+nq$ is a multiple of $n$, we have $H^{n+nq}(K(\ZZ,n);\ZZ)$ is infinite cyclic, generated by 
the cup product power $\iota_n^q$. The homotopy commutativity of \eqref{homotopy pullback} 
forces the underlying singular cohomology operation to be $\theta=\iota_n^q$  and therefore $\widehat{\theta}$ is a 
refinement of both the wedge power and cup product. It is known (see e.g. \cite{Bun1}) that the Deligne-Beilinson 
cup product is the unique refinement of these operations (up to homotopy) and $\widehat{\theta}={\rm DD}^q$.  

\medskip
For arbitrary $\lambda$, recall that the curvature map is surjective onto closed forms with integral periods. Hence, for $\lambda\not\in \ZZ$, the operation $\lambda \mu_n^q R$ is not in the image of $R$ and therefore $\tau$ does not admit a differential refinement. For $\lambda\in \ZZ$, $\lambda {\rm DD}^q$ defines a refinement and is again unique up to homotopy.

\vspace{0.1in}

\noindent (ii) \fbox{ $\tau=0$} 

\vspace{0.1in}

Let $\hat{\theta}$ be a refinement of $\tau=0$. Then $R\widehat{\theta}\simeq \tau R\simeq 0$. Since $\flat \BB^mU(1)$ is the fiber of $R$, $\widehat{\theta}$ must factor through the flat inclusion $j:\flat \BB^{m} U(1)\into \BB^mU(1)_{\nabla}$. Call the factorizing map $\phi^{\prime}$. By the homotopy commutativity of the diagram \eqref{homotopy pullback} and using Proposition \ref{refinement}, 
we also have that 
$$i\theta I\simeq i I\widehat{\theta}\simeq R\widehat{\theta}\simeq 0\;.$$
Hence, the image of $I$ must be killed by $i\theta$. But then homotopy commutativity implies that $i\theta$ must factorize through the point inclusion $\ast\to \flat_{\rm dR}\BB^mU(1)$. We see that we have a homotopy commutative diagram
 $$
 \xymatrix{
 \BB^{n+1}\ZZ \ar[d]^-{\theta}\ar[rr] && \ast \ar[d]
 \\
 \BB^{m+1}\ZZ\ar[rr]^-{i} && \flat_{\rm dR}\BB^{m+1}U(1)\;.
 }
 $$
Given the fiber sequence
 $$
 \flat \BB^{m} U(1)\overset{\beta}{\longrightarrow} \BB^{m+1}\ZZ\overset{i}{\longrightarrow} \flat_{\rm dR}\BB^{m+1}U(1)\;,
 $$
we see that the universal property gives a map $\phi:\BB^nU(1)\to \flat \BB^{m}U(1)$ such that $\theta\simeq \beta\phi$, 
where $\beta$ is the Bockstein map in the stacky diamond \eqref{stackdiamond}. Now, 
at the level of homotopy classes, we have
\bea
 \beta \phi I&=&\theta I   \qquad ~~  (\theta=\beta \phi~{\rm from~ above})
 \\
 &=& I \widehat{\theta}     \qquad  ~~({\rm Prop.}~ \ref{refinement})
 \\
 &=&Ij \phi^{\prime}  \qquad ({\rm factorization~ through~ flat~ inclusion})
 \\
 &=&\beta \phi^{\prime} \qquad~ ({\rm stacky ~diamond}~ \eqref{stackdiamond})\;,
 \eea
 which implies $\phi I-\phi^{\prime}$ is in the kernel of $\beta$. By exactness in the stacky diamond
 \eqref{stackdiamond}, this implies that $\phi I-\phi^{\prime}$ is in the image of
 \(
 {\rm exp}:\pi_0\map(\BB^nU(1)_{\nabla},\flat_{\rm dR}\BB^{m}U(1))\to \pi_0\map(\BB^nU(1)_{\nabla},\flat \BB^{m}U(1))\;.
 \label{exp map}
 \)
  If $m\neq kn$, for some $k>0$, then the group on the left 
  is zero by Lemma \ref{de rham ops}. 
  Hence, $\phi I=\phi^{\prime}$ and 
 $$
 j\phi I=j\phi^{\prime}=\widehat{\theta}\;.
 $$
 If $m=kn$, then the group on the left of \eqref{exp map}  is isomorphic to $\RR$, again by Lemma 
 \ref{de rham ops}. In this case, let us recall from Prop. \ref{Prop-disc} that we have an equivalence  
 $$
 \flat \BB^n U(1)\simeq \BB^nU(1)^{\delta}\simeq \flat \BB^nU(1)^{\delta}\;.
 $$
 Then, by cohesion, we have
 \begin{eqnarray}
 \pi_0\map(\BB^nU(1)_{\nabla},\flat \BB^{kn+1}U(1))&\simeq&  \pi_0\map({\bf \Pi}\BB^nU(1)_{\nabla},\BB^{kn}U(1)^{\delta})
 \quad~ ({\rm by~above~and~}\eqref{monad-adj})
\nonumber \\
  &\simeq & \pi_0\map(\BB^{n+1}\ZZ,\BB^{kn}U(1)^{\delta})
   \qquad \qquad ({\rm by~}\eqref{monad-adj})
\nonumber \\
 &\simeq & \pi_0\map(\BB^{n+1}\ZZ,\flat \BB^{kn}U(1))
  \qquad \qquad ({\rm by~above})\;.
 \label{kn map}
 \end{eqnarray}
Again, the isomorphism is provided by precomposing with $I: {\rm id} \to \Pi$ (eq. \eqref{eq I}). 
Now let $\varphi$ be such that $\phi I-\phi^{\prime}={\rm exp}(\varphi)$. By the above isomorphism \eqref{kn map}, 
it follows that there is
$$
\phi^{\prime\prime}\in \pi_0\map(\BB^{n+1}\ZZ,\flat \BB^{kn+1}U(1))\;,
$$
such that ${\rm exp}(\varphi)=\phi^{\prime\prime}I$. 
Hence, $\phi I-\phi^{\prime}={\rm exp}(\varphi)=\phi^{\prime\prime}I$, so that
 $\phi^{\prime}=(\phi-\phi^{\prime\prime})I$, since
$I$ is a linear operation. Applying $j$ to the latter equation gives the result.  
\endofproof

\section{Differential Steenrod operations}
\label{Sec Sq}

We now would like to apply and specialize the discussion in 
the previous section  to                                                                   
describe differential cohomology operations which refine the classical Steenrod squares
(see Sec. \ref{Sec Class}). That is, we seek  differential cohomology operations 
$\widehat{\theta}_i$ such that 
\(
\label{rho hat}
\rho_2 I \widehat{\theta}_k=Sq^k \rho_2 I\;.
\)
Here, $\rho_2: \Z \to \Z/2$ denotes the mod 2 reduction morphism.
Rephrasing this diagrammatically, we aim for a commutative diagram 
\(\label{refine steenrod}
\xymatrix{
\BB^{n} U(1)_{\nabla} \ar[d]^-{\widehat{\theta}_k} 
 \ar[r]^-{I} & \BB^{n+1}\ZZ \ar[r]^-{\rho_2} \ar@{-->}[d]& {\bf B}^{n+1} \ZZ/2 \ar[d]^-{Sq^k} 
 \\
\BB^{n+k} U(1)_{\nabla} \ar[r]^-{I} & \BB^{n+1+k}\ZZ\ar[r]^-{\rho_2} & {\bf B}^{n+1+k} \ZZ/2\;. 
}
\)
In the previous section, we saw that every differential cohomology operation refines a singular operation. 
Therefore, we can fill in the middle vertical arrow and ask for the entire diagram to commute up to homotopy. 

\medskip
However, as we saw in the Introduction, the homotopy commutativity of the right square is too much to ask in general. 
That is, not every $\ZZ/2$ operation admits an integral refinement. For example, the operations
$$
Sq^{2k}:H^n(-;\ZZ/2)\to H^{n+2k}(-;\ZZ/2)
$$
cannot have an integral refinement. Otherwise the operation would be in the image of the 
mod 2 reduction map and, by exactness, the Bockstein ${\beta}(Sq^{2m})$ 
(relating integral to mod 2 coefficients) would vanish. This is not the case, however, 
since the Adem relations imply
$$
0\neq Sq^{2k+1}=Sq^1Sq^{2m}=(\rho_2{\beta})Sq^{2k}\;.
$$ 
It therefore does not make sense to refine the even Steenrod squares. 
\begin{proposition}
The even Steenrod squares do not admit differential refinements. 
\label{Prop even}
\end{proposition}
Put another way, when refining mod 2 operations 
(or mod $p$ in general), one first needs an integral refinement. If such an integral lift exists, then one 
can ask for a differential refinement.

\medskip
Given the characterization theorem, Theorem \ref{characterization theorem}, established in the previous section, 
we can identify what these classes must be for odd Steenrod squares relatively easily. 
\begin{lemma}
\label{Prop mod p}
The odd integral Steenrod operations $Sq_{\ZZ}^{2k+1}:H^n(-;\ZZ) \to H^{n+2k+1}(-;\ZZ)$ factorizes uniquely as 
$$
\xymatrix{
\theta:H^n(-;\ZZ)\ar[r]^-{\rho_2} & 
H^n(-;\ZZ/2) \ar[r]^-{Sq^{2k}} &
H^{n+2k}(-;\ZZ/2) \ar[r]^-{\Gamma_2} & 
H^{n+2k}(-;U(1)) \ar[r]^-{\widetilde{\beta}} &
H^{n+2k+1}(-;\ZZ)
}\;,
$$
where $\Gamma_2$ is induced by the representation $\ZZ/2\into U(1)$ as the square roots of unity,
and $\widetilde{\beta}$ is the Bockstein corresponding to the exponential sequence 
$\Z \to \R \to U(1)$. 
\end{lemma}
\theproof
Recall that $Sq^{2k+1}_{\ZZ}$ is defined as the operation $\beta Sq^{2k}\rho_2$, where $\beta$ is the Bockstien corresponding to the sequence 
$$\ZZ\overset{\times 2}{\to} \ZZ \overset{\rho_2}{\to} \ZZ/2\;.$$
Now consider the morphism of short exact sequences
$$
\xymatrix{
0\ar[r] & \ZZ\ar[rr]^{\times 2\pi i } && \RR\ar[rr]^{{\rm exp}} && U(1)\ar[r] & 0
\\
0 \ar[r] & \ZZ\ar[rr]^{\times 2}\ar[u]^{\rm id} && \ZZ\ar[rr]^{\rho_2}\ar[u]^{\pi i\times } && \ZZ/2\ar[r]\ar[u]^{\Gamma_2} & 0~.
}
$$
This morphism induces a morphism of long fibration sequences involving the Bockstein homomorphisms 
$$
\xymatrix{
\hdots \ar[r] & B^{n-1}\ZZ\ar[rr]^{\times 2\pi i} && B^{n-1}\RR\ar[rr]^{{\rm exp}} && B^{n-1}U(1)\ar[rr]^-{\tilde{\beta}} 
&& B^n\ZZ\ar[r] & \hdots 
\\
\hdots \ar[r] & B^{n-1}\ZZ\ar[rr]^{\times 2}\ar[u]^{\rm id} && B^{n-1}\ZZ\ar[rr]^{\rho_2}\ar[u]^{\times \pi i } && B^{n-1}\ZZ/2\ar[rr]^-{\beta}\ar[u]^{\Gamma_2} && B^n\ZZ\ar[r]\ar[u]^{\rm id} & \hdots \;.
}
$$
The homotopy commutativity of the right square gives the desired factorization. Uniqueness follows from the definition of $Sq^{2k+1}_{\ZZ}$, along with the fact that every stable operation $\phi:B^{n-1}\ZZ/2\to B^{n-1}U(1)$ is induced by a representation $\Gamma_2:\ZZ/2\to U(1)$ (of which there is only 1).  
\endofproof

\medskip
As a corollary of Proposition \ref{Prop mod p} and the characterization theorem, Theorem \ref{characterization theorem}, 
we have the following.
\begin{corollary}
\label{Cor Sq odd}
Let $\widehat{\theta}_{2k+1}$ be a cohomology operation refining the odd integral Steenrod square
$Sq_{\ZZ}^{2k+1}$. Then we have
$$
\widehat{\theta}_{2k+1}=j\Gamma_2Sq^{2k}\rho_2 I\;,
$$
so that we can define the refinement as $\widehat{Sq}^{2k+1}:=\widehat{\theta}_{2k+1}$. Diagrammatically, we have 
$$
\hspace{-2mm}
\xymatrix{
\BB^nU(1)_{\nabla} \ar[rd]^{I} \ar[rrrrr]^-{\widehat{Sq}^{2k+1}} & & &  & & \BB^{n+2k+1}U(1)_{\nabla}\;.
\\
& \BB^{n+1}\ZZ\ar[rd]^{\rho_2} & &  & \flat \BB^{n+2k}U(1) \ar[ru]^{j} &
\\
&& \BB^{n+1}\ZZ/2 \ar[r]^-{Sq^{2k}}  & \BB^{n+2k}\ZZ/2 \ar[ru]^{\Gamma_2} & &
}$$
\end{corollary}
\theproof
Since $\widehat{\theta}$ refines $Sq_{\ZZ}^{2k+1}$, which takes values in torsion, we have $i(Sq^{2k+1}_{\ZZ})=0$. The homotopy commutativity of diagram \eqref{homotopy pullback} implies the corresponding operation on forms 
$\tau=0$. By theorem \ref{characterization theorem}, we must have that 
$$
\widehat{\theta}_{2k+1}=j\phi I\;,
$$
for some operation 
$$
\phi:\BB^{n+1}\ZZ\to \flat \BB^{n+2k} U(1)\;.
$$
Since $Sq_{\ZZ}^{2k+1} I=I\widehat{\theta}=\beta \phi I$, Proposition \ref{Prop mod p} implies that $\phi$ 
must be $\Gamma_2 Sq^{2k}\rho_2$.
\endofproof

\subsection{Relationship with the Deligne-Beilinson cup product}
\label{Sec DB}

In the previous section we established that the only differential refinement of the odd Steenrod squares is given by the operation
$\widehat{Sq}^{2k+1}:=j\Gamma_2 Sq^{2k} \rho_2 I$.
Hence, from the point of view of refinement our work is done. However classically, we know that the Steenrod squares are related to the homotopy commutativity of the cup product. One could ask whether or not the differential Steenrod squares are related to the homotopy commutativity of the Deligne-Beilinson cup product. This section can be viewed as a refinement of 
the second classical point of view on the Steenrod squares presented in Sec. \ref{Sec Class}.

\medskip
In fact, it is already known \cite{Go} that if $\hat{x}$ is a differential cohomology class of degree $2n+1$, then the Deligne-Beilinson square cup $\hat{x}^2$ is related to the image of the Steenrod square $Sq^{n-1}$ in differential cohomology via the map $j\Gamma_2$, introduced in the previous section.
At the end of the section, we generalize the result of Gomi \cite{Go}.

\medskip
Let $X$ be a manifold. As outlined in \cite{FSS1}, we have a cup product morphism in differential cohomology
$$
\xymatrix{
{X} \ar[r]^-{\Delta} & {X}\times {X} \ar[r]^-{\hat{x}\times \hat{x}} &
\BB^{n}U(1)_{\nabla}\times \BB^{n}U(1)_{\nabla} \ar[r]^-{\cup_{\rm DB}} & \BB^{2n+1}U(1)_{\nabla}
}\;.
$$
As in the classical case, the Deligne-Beilinson cup product is not strictly graded commutative, but is graded commutative up to homotopy.  That is, we have a homotopy commutative diagram in stacks
$$
\xymatrix@C=0em @R=1em{
 & {X}\times {X}\ar[rrrdd]^-{\hat{x}\times \hat{x}}\ar[dddd] &&&&&&&&
\\
(x,y) \ar@{|->}[dd] &&&&&&&&&
 \\
 &&&&   \BB^{n}U(1)_{\nabla}\times \BB^{n}U(1)_{\nabla}  \ar[rrrrr]^-{\cup_{\rm DB}} &&&&& \BB^{2n+1}U(1)_{\nabla}
 \\
  (y,x) &&&&&&&&&
 \\
 & {X}\times {X}\ar[rrruu]_-{\hat{x}\times \hat{x}} &&&&&&&&&
 }
 $$
If we choose homotopies and higher coherence homotopies filling the diagram, we can equivalently express this by saying that $\BB^{2n+1}U(1)_{\nabla}$ is an $(\infty,1)$-cocone over the the diagram given by the $\ZZ/2$-action (call it $\psi$)
 on $X\times X$ via the above transposition map. 
If we take the colimit over this $\ZZ/2$-action, then the universal property of the colimit will ensure that there is a map (unique up to homotopy) from this colimit to $\BB^{2n+1}U(1)_{\nabla}$. That is, 
we have the following. 

\begin{lemma} 
The colimit of the $\Z_2$-action $\psi$ (described above)  sits a homotopy commuting diagram
\begin{center}
\begin{tikzpicture}
\matrix (m) [matrix of math nodes, column sep=3em, row sep=3em]{
 {\rm hocolim}(\psi)&
\\
{X} \times {X} & \BB^{2n+1}U(1)_{\nabla}\;.
\\
};

\path[->] 
(m-2-1) edge (m-1-1)
(m-2-1) edge node[above] {\footnotesize $\cup$} (m-2-2)
(m-1-1) edge[dashed] node[above] {\footnotesize $\hat{\lambda}$} (m-2-2);
\end{tikzpicture}
\end{center}
\end{lemma}

\begin{remark}
The colimit here serves to extract the homotopies involved in the $\ZZ/2$-action. The map $\hat{\lambda}$ attaches homotopies involved with the cup product to these homotopies. 
\end{remark}

Although there may be several ways to compute this colimit, we will make use of the cohesive structure on smooth stacks to perform the calculation.
\begin{proposition}
Let ${\bf Y}$ be a stack equipped with an action of $\ZZ/2$, that is, a functor $\psi:\ZZ/2 \to Sh_{\infty}(\cartsp)$ 
sending the unique object $\ast\in \ZZ/2$ to the stack ${\bf Y}$. The colimit over this functor is computed as
$$
{\rm hocolim}(\psi)\simeq {\bf E} \ZZ/2\times_{\psi} {\bf Y}\;,
$$
where ${\bf E} \ZZ/2={\rm disc}(S^{\infty})$ is the discrete universal principal $\ZZ/2$-bundle over the discrete stack $\BB \ZZ/2={\rm disc}(\RR P^{\infty})$.
\end{proposition}
\theproof
Since  the prestack category $[\cartsp,\sset]$ is combinatorial and simplicial, the homotopy colimit in prestacks is presented by the local homotopy colimit 
$$
{\rm hocolim}_{\rm local}(\psi)=\int^{\ast \in \ZZ/2}{\cal N}((\ZZ/2)/\ast)\odot \psi(\ast)\;.
$$
Here, ${\cal N}$ denotes the nerve while $\odot$ denotes the tensoring of a stack and a simplicial set. To compute the
 right hand side, we observe that the tensoring of a prestack ${\bf Y}$ and a simplicial set $X$ is provided by taking 
 the product with the constant stack 
$$
X\odot {\bf Y}:= {\rm const}(X)\times {\bf Y}\;.
$$
Then the coend is computed as
\bea 
\int^{\ast \in \ZZ/2}{\cal N}((\ZZ/2)/\ast)\odot \psi (\ast) &= & \int^{\ast \in \ZZ/2} E\ZZ/2 \odot {\bf Y}
\\
&=& \int^{\ast \in \ZZ/2} {\rm const}(E\ZZ/2) \times {\bf Y}
\\
&=& {\rm coeq}\Big\{\xymatrix{ {\rm const}(E\ZZ/2) \times {\bf Y}\ar@<.25em>[r]^{{\rm id}}\ar@<-.25em>_{\psi}[r] &  {\rm disc}(E\ZZ/2) \times {\bf Y}}\Big\}
\\
&=& {\rm const}(E \ZZ/2)\times_{\psi} {\bf Y}\;.
\eea
The homotopy colimit was computed in prestacks. Since the stackification functor is a left $\infty$-adjoint, it preserves homotopy colimits and we need only compute the stackification of the prestack ${\rm const}(E \ZZ/2)\times_{\psi} {\bf Y}$. 
Since ${\bf Y}$ was assumed to be a stack, this is ${\rm disc}(E\ZZ/2)\times_{\psi}{\bf Y}$, as claimed.
\endofproof

\begin{corollary} For the trivial action $\psi$ on a stack ${\bf Y}$, we have
$$
{\rm hocolim}(\psi)\simeq {\rm disc}(B\ZZ/2)\times {\bf Y}\simeq 
\BB\ZZ/2\times {\bf Y}\;.
$$
\end{corollary}

Returning to our discussion, we can now unravel the homotopies contained in the $\ZZ/2$-action.

\begin{proposition} 
\label{Prop cup h}
The stacky cup product map
${X} \to {X} \times {X} \to \BB^{2n+1}U(1)_{\nabla}$
can be extended to a map $\hat{\lambda}$ making the diagram
\begin{center}\label{homotopy diagram}
\begin{tikzpicture}
\matrix (m) [matrix of math nodes, column sep=3em, row sep=3em]{
{X}\times \BB\ZZ/2 & { X} \times {X} \times_{\ZZ/2} {\bf E}\ZZ/2&
\\
{X} & {X} \times {X} & \BB^{2n+1}U(1)_{\nabla}
\\
};

\path[->] 
(m-1-1) edge node[above] {\footnotesize $Q(\Delta)$} (m-1-2)
(m-2-1) edge (m-1-1)
(m-2-2) edge (m-1-2)
(m-2-1) edge node[above] {\footnotesize $\Delta$} (m-2-2)
(m-2-2) edge node[above] {\footnotesize $\cup$} (m-2-3)
(m-1-2) edge[dashed] node[above] {\footnotesize $\hat{\lambda}$} (m-2-3);
\end{tikzpicture}
\end{center}
commute up to homotopy. Moreover, given choices of homotopies and higher homotopies filling the diagram, $\hat{\lambda}$ is uniquely determined up to homotopy. Here, the two vertical maps are canonical sections of the projection $p_{X}:{X}\times \BB/\ZZ/2 \to {X}$ and $Q(\Delta)$ is yet to be determined.
\end{proposition}
\theproof
Equip ${X}$ with the trivial $\ZZ/2$-action, and equip ${X}\times {X}$ with the action given by transposing the two factors. Then the diagonal 
$$
\Delta:{X}\to {X}\times {X}
$$
defines a natural transformation of $\ZZ/2$-actions, and hence induces a map $Q(\Delta)$ on the corresponding homotopy colimits. Moreover, by the homotopy commutativity of the cup product, the map
$$
\xymatrix{
{X}\times {X} \ar[r]^-{\hat{x}\times \hat{x}} &
\BB^{n}U(1)_{\nabla}\times \BB^{n}U(1)_{\nabla} \ar[r]^-{\cup_{\rm DB}} & \BB^{2n+1}U(1)_{\nabla}
}
$$
commutes (up to homotopy) with the $\ZZ/2$-action. Given a choice of homotopies and higher homotopies, the universal property for $(\infty,1)$-colimits produces a map $\hat{\lambda}$, defined uniquely up to homotopy, making the diagram commute.
\endofproof

\medskip
To extract the Steenrod squares from this diagram, we will need to choose homotopies filling the diagram and study the composite map $\hat{\lambda}Q(\Delta)$. This is analogous to the classical case, where one produces such a diagram and then used the K\"unneth formula to compute the degree $2n$ cohomology of $X\times \RR P^{\infty}$. The coefficients are then defined to be the Steenrod squares. 

\begin{remark}
It is interesting to note that our method seems conceptually much simpler than the classical construction. However, we emphasize the fact that the explicit construction of the higher coherence homotopies would be just as complicated as in the classical case. Fortunately, we will be able to use the classical construction to our advantage for the choice of homotopies. 
\end{remark}

\medskip
As indicated, we will need to make use of a K\"unneth-type theorem for differential cohomology. Although it is likely that such a theorem follows from a more general theorem for sheaf hypercohomology, this particular case does not require such machinery and we can prove the claim directly \footnote{In the journal version of this article, the formulation of the K\"unneth theorem is false. The correct hypothesis on $Y$ is that it is a discrete simplicial presheaf of finite type. We have corrected this error in the present version.} .
\begin{proposition}(K\"unneth decomposition for differential cohomology)
Let $X$ be a smooth manifold and let $Y={\rm disc}(Y')$, where $Y'$ is a simplicial set of finite type (only finitely many simplices in each simplicial degree). Then we have a natural short exact sequence
$$
\begin{tikzpicture}[descr/.style={fill=white,inner sep=1.5pt}]
        \matrix (m) [
            matrix of math nodes,
            row sep=1em,
            column sep=0em,
            text height=1.5ex, text depth=0.25ex
        ]
        { 0 &   \widehat{H}^n(X;\ZZ) \oplus \bigoplus_{i=1}^n H^{n-i-1}(X;U(1))\otimes H^{i}(Y;\ZZ)  & \widehat{H}^{n}(X\times Y;\ZZ)   &\\ 
            & {\rm Tor}(\widehat{H}^n(X;\ZZ),H^1(Y;\ZZ))\oplus \bigoplus_{i=1}^{n} {\rm Tor}\left(H^{n-i-1}(X;U(1)), H^{i+1}(Y;\ZZ) \right)  & 0~. \\
        };

        \path[overlay,->, font=\scriptsize,>=latex]
        (m-1-1) edge (m-1-2)
        (m-1-2) edge (m-1-3)
        (m-1-3.east) edge[out=350,in=170] (m-2-2.west)
        (m-2-2) edge (m-2-3);
\end{tikzpicture}
$$
Moreover, the sequence splits (but not naturally).
\label{Prop Kunneth}
\end{proposition}
\theproof
Let $\{U_i\}$ be a good open covers of $X$ and let $C(\{U_i\})$ be the {\v C}ech nerve of the cover. Since $Y={\rm disc}(Y')$ for some simplicial set $Y'$, it is cofibrant in the projective model structure on simplicial presheaves. Since this model category is cartesian, the product $C(\{U_i)\})\times Y$ is a projective resolution of $X\times Y$. 

The total complex of the {\v C}ech-Deligne double complex  of $X\times Y$ is, by definition, the hom in unbounded chain complexes
$$
{\rm tot}\left(C_{\bullet\bullet}(C(\{U_i\})\times Y);\ZZ^{\infty}_{\cal D}(n+1))\right):=\hom_{\ch}\left(C_{\bullet}(C(\{U_i\}) \times Y)), \ZZ^{\infty}_{\cal D}(n+1)\right)\;.
$$
where 
$C_{\bullet}:{\rm PSh}_{\infty}(\cartsp) \to {\rm PSh}(\cartsp;\sab)\to {\rm PSh}_{\infty}(\cartsp;\mathscr{C}{\rm h}_{+})$ 
is the (prolongation to presheaves) of the composition of the free functor and the Moore complex functor. Now the Eilenberg-Zilber map gives a chain homotopy equivalence
$$\nabla: C_{\bullet}(C(\{U_i\}))\otimes C_{\bullet}(Y)
\overset{\simeq}{\longrightarrow} C_{\bullet}(C(\{U_i\})\times Y)\;.$$
Taking the hom in chain complexes to the Deligne complex 
it follows that we have an equivalence
$$
\nabla:\hom_{\ch}\left(C_{\bullet}(C(\{U_i\})\times Y),\ZZ^{\infty}_{\cal D}(n+1))\right)\overset{\simeq}{\to} \hom_{\ch}\left(C_{\bullet}(\{U_i\})\otimes C_{\bullet}(Y),\ZZ^{\infty}_{\cal D}(n+1))\right)\;.
$$
Since $Y'$ is finitely generated in each simplicial degree, the appropriate finiteness assumption is satisfied and the canonical map
$$ \hom_{\ch}\left(C_{\bullet}(\{U_i\}),\ZZ^{\infty}_{\cal D}(n+1))\right)\otimes \hom_{\ch}\left(C_{\bullet}(Y),\ZZ)\right)\to \hom_{\ch}\left(C_{\bullet}(\{U_i\})\otimes C_{\bullet}(Y),\ZZ^{\infty}_{\cal D}(n+1))\right)$$
is a quasi isomorphism.

Write
$$C^j(Y;\ZZ)=\hom_{\ch}\left(C_{\bullet}(Y),\ZZ)\right)_{-j}$$
for the {\v C}ech cochain complex of $Y$. For the total {\v C}ech Deligne complex, we have
$$
{\rm tot}_{n-j}\left(C^{\bullet\bullet}(\{U_i\};\ZZ^{\infty}_{\cal D}(n+1))\right)=
\hom_{\ch}\left(C_{\bullet}(\{U_i\}),\ZZ^{\infty}_{\cal D}(n+1))\right)_{j}\;.
$$
Then 
\begin{align} 
\hspace{-.6cm}
\phantom{i+j+k}
& \begin{aligned}
\mathllap{
 \bigoplus_{i+j=0}} & H_i\left(\hom_{\ch}\left(C_{\bullet}(\{U_i\}),\ZZ^{\infty}_{\cal D}(n+1)\right)\right)\otimes H_j\left(\hom_{\ch}(C_{\bullet}(Y),\ZZ)\right) 
 \\
 & \qquad \hspace{6.5cm}
= \bigoplus_{i+j=0} {\bf H}^{n-i-1}(X;\ZZ^{\infty}_{\cal D}(n+1))\otimes H^{-j}(Y;\ZZ)
\\
 & \qquad  \hspace{6.5cm}
= \widehat{H}^n(X;\ZZ) \oplus \bigoplus_{i=1}^n H^{n-i-1}(X;U(1))\otimes H^{i}(Y;\ZZ)
\end{aligned}
\nonumber
\end{align}
and
 \begin{align}
 \hspace{.2cm}
   \phantom{i + j + k}
 &\begin{aligned}
    \mathllap{
 \bigoplus_{i+j=-1}{\rm Tor}} & \left(H_i\left(\hom_{\ch}(C_{\bullet}(\{U_i\}),\ZZ^{\infty}_{\cal D}(n+1))\right), H_j\left(\hom_{\ch}(C_{\bullet}(Y),\ZZ)\right)\right) 
 \\
& \qquad \hspace{2cm} 
= \bigoplus_{i+j=-1} {\rm Tor}\left({\bf H}^{n-i}(X;\ZZ^{\infty}_{\cal D}(n+1)), H^{-j}(Y;\ZZ)\right)
 \\
 & \qquad \hspace{2cm} 
 = {\rm Tor}(\widehat{H}^n(X;\ZZ),H^1(Y;\ZZ))\oplus \bigoplus_{i=1}^{n} {\rm Tor}\left(H^{n-i-1}(X;U(1)), H^{i+1}(Y;\ZZ) \right)\;.
 \end{aligned}
 \nonumber
 \end{align} 
 Now the result follows from the usual K\"unneth formula.
\endofproof

We now adapt the general K\"unneth decomposition to the case directly related to Steenrod squares.  

\begin{proposition}\label{cup homotopy}
For a smooth manifold $X$, we have
$$
\widehat{H}^{2n}(X\times \BB\ZZ/2;\ZZ)\simeq \widehat{H}^{2n}(X;\ZZ)\oplus \bigoplus_{j<2n\ {\rm even}}T^j_2\;,
$$
where $T^i_2$ is the 2-torsion subgroup of $\widehat{H}^i(X;\ZZ)$.
\end{proposition}
\theproof
Applying Proposition \ref{Prop Kunneth} with $Y=\mathbf{B}\ZZ/2$, we have
\begin{align} 
\hspace{-.8cm}
\phantom{i+j+k}
& \begin{aligned}
\mathllap{\widehat{H}^{2n}}&({X}\times \BB\ZZ/2;\ZZ)
\cong \widehat{H}^{2n}(X;\ZZ) \oplus \bigoplus_{i=1}^n H^{2n-i-1}(X;U(1))\otimes 
H^{i}(\BB\ZZ/2;\ZZ)\oplus
\\
 & \qquad \hspace{2cm} 
 \oplus  {\rm Tor}(\widehat{H}^{2n}(X;\ZZ),H^1(\BB\ZZ/2;\ZZ))\oplus \bigoplus_{i=1}^{2n} {\rm Tor}\left(H^{2n-i-1}(X;U(1)), 
 H^{i+1}(\BB\ZZ/2;\ZZ) \right)\;.
 \end{aligned}
 \nonumber
 \end{align} 
It remains to identify these groups. We first note that we have an isomorphism
\(
\label{u1 trivial} H^i(X,U(1))\otimes \ZZ/2\simeq 0
\)
for all $i$. To see this, consider the short exact sequence 
$\ZZ \overset{\times 2}{\longrightarrow} \ZZ \to \ZZ/2$.
Tensoring with $H^i(X;U(1))$ leads to the sequence
$$
H^i(X,U(1))\overset{\times 2}{\longrightarrow} H^i(X;U(1))\to H^i(X;U(1))\otimes \ZZ/2\to 0\;.
$$
But since we have $U(1)$ coefficients, the map $\times 2$ is surjective and the first isomorphism theorem 
confirms the claim \eqref{u1 trivial}. 

Now recall the integral cohomology groups of the classifying space $K(\Z_2, 1)=\RR P^{\infty}$ 
$$H^i(\RR P^{\infty};\ZZ)=\left\{\begin{array}{ccl}
\ZZ && i=0,
\\
\ZZ/2 && i\ {\rm even}\neq 0,
\\
0 && {\rm otherwise}.
\end{array}
\right.$$
Combining this with equation \eqref{u1 trivial} gives 
$$\widehat{H}^{2n}(X;\ZZ) \oplus \bigoplus_{i=1}^n H^{2n-i-1}(X;U(1))\otimes H^{i}(\BB\ZZ/2;\ZZ)\simeq \widehat{H}^{2n}(X;\ZZ)\;.$$
Then the Tor groups are easily computed
\begin{align} 
\hspace{-.6cm}
\phantom{i+j+k}
& \begin{aligned}
\mathllap{
 {\rm Tor}} & (\widehat{H}^{2n}(X;\ZZ),H^1(\BB\ZZ/2;\ZZ))\oplus \bigoplus_{i=1}^{2n} {\rm Tor}\left(H^{2n-i-1}(X;U(1)), H^{i+1}(\BB\ZZ/2;\ZZ) \right)  
 \\
 & \qquad \hspace{7.5cm}
\simeq  \bigoplus_{1\leq i\leq 2n,\ {\rm odd}} {\rm Tor}\left(H^{2n-i-1}(X;U(1)), \ZZ/2) \right) 
\\
 & \qquad \hspace{7.5cm}
 \simeq \bigoplus_{j<2n\ {\rm even}}T^j_2\;.
\end{aligned}
\nonumber
\end{align}

\vspace{-5mm}
\endofproof

In fact, we can be a bit more precise about what the torsion groups $T^j_2$ actually look like. 
To that end, let us start by recalling the following.  
\begin{lemma} 
As a ring, the integral cohomology of $\mathbf{B}\ZZ/2$ takes the form
$$
H^*(\mathbf{B}\ZZ/2;\ZZ)\cong H^{\ast}(\RR P^{\infty};\ZZ)\simeq \ZZ[x]/\langle 2x\rangle\;,
$$
where $x$ is an integral lift of $w_1$.
\end{lemma}

The statement of the lemma is classical and can be established in various ways. For example, one can view
$\R P^\infty$ as a Grassmannian 
$G_{2m+1}(\R^\infty)$ with $m=0$, for which 
the Bockstein exact sequence 
$$
\xymatrix{
\hdots \ar[r] & H^j(-; \Z) \ar[r]^-{\times 2} &
H^j(-; \Z) 
\ar[r]^-{\rho_2}
&
H^j(-; \Z/2)
\ar[r]^-\beta
&
H^{j+1}(-; \Z) \ar[r] & \hdots
}
$$
implies that the integral cohomology $H^*(G_{2m+1}(\R^\infty); \Z)$ 
splits additively as the direct sum of a polynomial $\Z[p_1, \cdots, p_m]$
and the image of $\beta$ (see \cite{MiS} Problem 15-C). For $m=0$ this then gives 
that the integral cohomology of $\RP^\infty$ is the image of $\beta$. 
This can also be deduced via chain complexes 
(see \cite{Hat} p. 222). A third way is to consider the Gysin sequence 
for integral cohomology corresponding to the circle bundle 
$S^1 \to S^\infty \overset{\pi}{\longrightarrow} \RP^\infty$, 
that is 
$$
\xymatrix{
\hdots \ar[r] &
 H^n(S^\infty; \Z) \ar[r]^-{\pi_*} &
H^{n-1}(\RP^\infty; \Z) \ar[r]^-{\cup e} &
H^{n+1}(\RP^\infty; \Z) \ar[r]^-{\pi^*} &
H^{n+1}(S^\infty;\Z) \ar[r] & \hdots\;,
}
$$
where $\pi^*$ is pullback, $\pi_*$ is pushforward, and $e$ is the Euler class
of the circle bundle. The latter gives an isomorphism between all even degree cohomology
groups of $\RP^\infty$. Then the Euler class $e$, being 2-torsion, gives the desired 
result.

\medskip
Going back to the 2-torsion subgroup $T_2^j$ , for even $j$, the torsion pairing is given explicitly as
$$
{\rm Tor}\left(H^{2n-i-1}(X;U(1)), H^{i+1}(\mathbf{B}\ZZ/2;\ZZ) \right)\simeq {\rm Tor}
\left(H^{2n-j}(X;U(1)), \ZZ/2\langle x^j\rangle \right)\;.
$$
Now the sequence
$$
0 \to \ZZ\langle x^j\rangle \overset{\times 2}{\longrightarrow} \ZZ\langle x^j\rangle \to \ZZ/2\langle x^j\rangle \to 0
$$
is a projective resolution. Therefore, the torsion pairing is given as the kernel of the map
$$\times 2:H^{2n-j}(X;U(1))\otimes \ZZ\langle x^j\rangle\to H^{2n-j}(X;U(1))\otimes \ZZ\langle x^j\rangle$$
which is spanned by elements of the form $y\otimes x^j$, with $y$ a 2-torsion element in $H^{2n-j}(X;U(1))$. 

\medskip
Finally, let us return to our analysis of the map
$$
\hat{\lambda}Q(\Delta):{X}\times \BB \ZZ/2 \to \BB^{2n}U(1)_{\nabla}\;,
$$
defined in the diagram of Proposition \ref{Prop cup h}. By Proposition \ref{cup homotopy} and the above
 discussion, we see that we 
 can expand the class of the map $\hat{\lambda}Q(\Delta)$ 
  as a homogeneous polynomial
\(\label{polynomial}
[\hat{\lambda}Q(\Delta)]=\widehat{s}_n^{~2n}+  \widehat{s}_n^{~2n-2}\otimes x+\hdots + \widehat{s}_n^{~2}\otimes x^{2n-2}+\widehat{s}_n^{~0}\otimes x^{2n}\;,
\)
where each $\widehat{s}^{~2k}$ represents a differential cohomology operation. Although we know the general form that the map $\hat{\lambda}Q(\Delta)$ takes, the homotopy class of this map still depends on an explicit choice of homotopies and higher homotopies. For now, we leave these choices implicit and return to this point later.

\begin{definition} 
Define the operations/classes $\widehat{s}_n^{~2k}:\BB^nU(1)_{\nabla} \to \BB^kU(1)_{\nabla}$ by the expansion \eqref{polynomial}.
\end{definition} 

\begin{remark} {\bf (i)} 
Notice that naturality of the cup product implies that the classes $\widehat{s}_n^{~2k}$ are natural with respect to pullback and hence define differential cohomology operations. 

\medskip
\noindent {\bf (ii)} Notice that for $k<n$, the operation
$\widehat{s}_n^{~2k}$ 
must represent a trivial class. Indeed, since $\widehat{s}_n^{~2k}$ has image in 2-torsion the curvature
vanishes
$$
2R(\widehat{s}_n^{~2k}) =R(2\widehat{s}_n^{~2k})=R(0)=0\;,$$
indicating that the class $\widehat{s}_n^{~2k}$ takes values in flat bundles. It follows that the map $\widehat{s}_n^{~2k}$ factorizes through the stack $\BB^kU(1)^{\delta}$, representing cohomology with $U(1)$-coefficients. Since there are no degree-decreasing cohomology operations with $U(1)$-coefficients, the class of $\widehat{s}_n^{~2k}$ must be trivial in this case.
\end{remark}

 \medskip
It remains to identify the classes $\widehat{s}_n^{~2k}$ for $k>n$. We start with the top class. 
\begin{proposition}
The class $\widehat{s}_n^{~2n}(\hat{x})$ defined by the polynomial expression 
\eqref{polynomial} can be identified with the cup product $\hat{x}\cup\hat{x}$. 
\end{proposition}
\theproof
By the homotopy commutativity of the diagram \eqref{homotopy diagram}, 
the pullback of the class $[\hat{\lambda}Q(\Delta)]$ by the canonical section ${X}\to {X}\times \BB\ZZ/2$ is the cup product. This pullback simply restricts a class in $\widehat{H}^{2n}(X\times \BB/\ZZ/2;\ZZ)$ to $X$. From the polynomial expansion of $[\hat{\lambda}Q(\Delta)]$, it is apparent that this class is $\widehat{s}_n^{~2n}(\hat{x})$.
\endofproof

Recalling that the classes $\widehat{s}_n^{~2l}$ were left undetermined, we have the following.
\begin{proposition}
\label{Sq odd s}
There is a choice of homotopy commutative diagram \eqref{homotopy diagram} such that
$$
\widehat{s}_{n}^{~2k}=\widehat{Sq}^{2k+1}=j\Gamma_2Sq^{2k}\rho_2I\;.
$$
Moreover, when $n$ is odd, these homotopies uniquely refine the homotopies involved in the classical case.
\end{proposition}
\theproof
Since the map $\hat{\lambda}Q(\Delta)$ was left ambiguous, we simply define $\widehat{s}^{~2k}_n$ as needed. The homotopy class of this map then determines a homotopy commutative diagram \eqref{homotopy diagram}.

\medskip
To see how these relate to the homotopies involved in the classical case. Observe that since the DB cup product refines the singular cup product, we have a homotopy commutative diagram
$$\hspace{-.2cm}
\begin{tikzpicture}\label{diagram1} 
\matrix (m) [matrix of math nodes, column sep=3em, row sep=3em]{
{X}\times \BB\ZZ/2 & {X} \times {X} \times_{\ZZ/2} {\bf E}\ZZ/2 & &
\\
{X} & {X} \times {X} & \BB^{n-1}U(1)_{\rm conn}\times \BB^{n-1}U(1)_{\nabla} & \BB^{2n-1}U(1)_{\nabla}
\\
 & {X}\times {X} & \BB^{n}\ZZ/2\times \BB^n\ZZ/2 & \BB^{2n}\ZZ/2\;.
 \\
};

\path[->] 
(m-1-1) edge node[above] {\footnotesize $Q(\Delta)$} (m-1-2)
(m-2-1) edge (m-1-1)
(m-2-2) edge (m-1-2)
(m-2-1) edge node[above] {\footnotesize $\Delta$} (m-2-2)
(m-2-2) edge node[above] {\footnotesize $\hat{x}\times\hat{x}$} (m-2-3)
(m-2-3) edge node[above] {\footnotesize $\cup_{\rm DB}$} (m-2-4)
(m-1-2) edge[dashed, out=0, in=150] node[above] {\footnotesize $\hat{\lambda}$} (m-2-4)
(m-2-1) edge (m-3-2)
(m-2-2) edge (m-3-2)
(m-3-2) edge node[below] {\footnotesize $\rho_2\circ I(\hat{x})\times \rho_2\circ I(\hat{x})$} (m-3-3)
(m-2-3) edge node[right] {\footnotesize $\rho_2\circ I\times \rho_2\circ I$} (m-3-3)
(m-2-4) edge node[right] {\footnotesize $\rho_2\circ I$} (m-3-4)
(m-3-3) edge node[above] {\footnotesize $\cup$} (m-3-4);
\end{tikzpicture}
$$
Using the polynomial expansion of $[\hat{\lambda}Q(\Delta)]$, we can write the homotopy class of the 
upper top-left to bottom-right composite as
\begin{alignat}{2}
\rho_2I[\hat{\lambda}Q(\Delta)] &= [\rho_2I\hat{\lambda}Q(\Delta)]
\nonumber \\
&=\rho_2I\big(\widehat{s}_n^{~2n}+\widehat{s}_n^{~2n-2}\otimes x+\hdots + \widehat{s}_n^{~2}\otimes x^{n-1}
+\widehat{s}_n^{~0}\otimes x^{n}\big)
\nonumber \\
&= \rho_2I\left(\widehat{s}_n^{~2n}+\widehat{s}_n^{~2n-2}\otimes x+\hdots + \widehat{s}_n^{~n+1}\otimes x^{(n-1)/2}+
\widehat{s}_n^{~n}\otimes x^{n/2}+\hdots  \right)
\nonumber \\
&= \left(\rho_2I \widehat{s}_n^{~2n}\right)+\left(\rho_2I\widehat{s}_n^{~2n-2}\right)\otimes w_1^2+\hdots +
 \left(\rho_2I\widehat{s}_n^{~n+2}\right)\otimes w_1^{n-2}+\left(\rho_2I\widehat{s}_n^{~n}\right)\otimes w_1^{n}\;.
\nonumber
\end{alignat} 
Now the using the classical construction of the Steenrod squares discussed in we recall that there is a map from the top-left corner to the bottom-right given by
$$
[Q(\Delta)\lambda]=Sq^{n}\rho_2I+ Sq^{n-1}\rho_2I\otimes w_1+\hdots +Sq^1\rho_2I\otimes w_1^{n-1}+w_1^{n}\;.
$$
We would like to compare this polynomial with the previous one to identify the coefficients. Unfortunately, the map $Q(\Delta)\lambda$ can not be homotopic to $\rho_2I\hat{\lambda}Q(\Delta)$. This is immediately clear from the fact that $Sq^{k}$ is not in the image of the mod 2-reduction for even $k$.  This also reflects the fact that we cannot choose homotopies and higher homotopies filling the top diagram which are mapped to the right homotopies in the classical, outer diagram. However, we can split the map $Q(\Delta)\lambda$ into two parts depending on the parity of the exponent of $Sq^k$. That is, we define
$$
\lambda_0=\sum_{k\ \text{even},\ k\leq n}Sq^k\otimes w_1^{n-k} \qquad
{\rm and} \qquad
\lambda_1=\sum_{k\ \text{odd},\ k\leq n}Sq^k\otimes w_1^{n-k}\;.
$$
Now recall that for an odd Steenrod square $Sq^{2k+1}$, we have $Sq^{2k+1}=Sq^1Sq^{2k}$. Since $Sq^1Sq^1=0$, we have that $Sq^{2k+1}$ is in the kernel of $Sq^1$. The equation $Sq^1=\rho_2\beta$ relating $Sq^1$ to the Bockstein $\beta$ for the mod 2 reduction $\ZZ\to \ZZ/2$ implies that $Sq^{2k+1}$ must be in the image of the mod 2 reduction $\rho_2$. Moreover, since $n$ is odd, $w_1^{n-k}$ is an even power when $k$ is odd and $w_1^{n-k}$ is in the image of the mod 2-reduction. Factorizing the map $\lambda_1$ through the mod 2-reduction $\rho_2$ and integration map $I$, we can write 
$$
\lambda_1=\rho_2 I \alpha\;.
$$
 By the universal property, $\alpha$ determines a homotopy commutative diagram and, by expression \eqref{polynomial},
  is a map of the form 
 $$
 \alpha=\widehat{s}_n^{~2n}+\widehat{s}_n^{~2n-2}\otimes x+\hdots + \widehat{s}_n^{~2}\otimes x^{n-1}
 +\widehat{s}_n^{~0}\otimes x^{n}\;.
 $$
Setting $\hat{\lambda}Q(\Delta):=\alpha$, we have
$$
\rho_2 I [\hat{\lambda}Q(\Delta)]=[\lambda_1]\;.
$$
Comparing coefficients, we see that
$$
Sq^{2k+1}\rho_2I=Sq^1 Sq^{2k}\rho_2I=\rho_2 \beta Sq^{2k}\rho_2I=\rho_2 Ij Sq^{2k}\rho_2I=\rho_2 I \widehat{s}_n^{~2k}\;.
$$
Since both $Sq^{2k}\rho_2I$ and $\widehat{s}_n^{~2k}$ take values in 2-torsion, we must have $jSq^{2k}\rho_2 I=\widehat{s}_n^{~2k}$. Alternatively, since $\widehat{s}_n^{~2k}$ refines $Sq^{2k}$, we could use Theorem \ref{characterization theorem} to conclude that $\widehat{s}_n^{~2k}$ has the desired form. 
\endofproof

\medskip
In \cite{Go}, it was observed that for an odd degree differential cohomology class $\hat{x}$, the Deligne-Beilinson square is given by the inclusion of the $(n-1)$st Steenrod square $Sq^{n-1}\rho_2I(\hat{x})$ into differential cohomology via the representation of $\Gamma_2:\ZZ/2\into U(1)$ as the primitive square roots of unity. We provide another proof of this fact to highlight the power of the stacky perspective.

\begin{proposition}
\label{Prop trapezoid} 
For each $n>0$, the Steenrod square $Sq^{2n}$  fits into a homotopy commutativity of the diagram:
$$
\hspace{-1mm}
\xymatrix{
\BB^{2n}U(1)_{\nabla} \ar[rd]^{I} \ar[rrrrr]^{\cup_{\rm DB}^2}  & & & & & \BB^{4n+1}U(1)_{\nabla}
\\
& \BB^{2n+1}\ZZ\ar[rd]^{\rho_2} & &  & \flat \BB^{4n+1}U(1) \ar[ru]^{j} &
\\
&& \BB^{2n+1}\ZZ/2 \ar[r]^-{Sq^{2n}}  & \BB^{4n+1}\ZZ/2 \ar[ru]^{\Gamma_2} & &
}\;.$$
\end{proposition}
\theproof
Let ${\rm DD}_{2n+1}$ denote the Dixmier-Douady class in degree $2n+1$. Since this class has odd degree, the homotopy commutativity of the Deligne-Beilinson (DB) cup product implies that its square is 2-torsion, $2{\rm DD}^2=0$.
Consequently,  the curvature obeys
$$
R(2{\rm DD}^2)=2R({\rm DD})^2=0\;.
$$
This, in turn implies $R({\rm DD})^2=0$. 
It follows that the square cup factorizes as
$$
{\rm DD}^2:\BB^{2n}U(1)_{\nabla}\to \flat \BB^{4n+1}U(1) \to \BB^{4n+1}U(1)_{\nabla}\;.
$$
Since ${\rm DD}^2$ is 2-torsion, it is killed by the $\times 2$ map. Given the Bockstein sequence associated to
$$
0\to \ZZ/2\to U(1) \overset{\times 2}{\longrightarrow} U(1)\to 0\;,
$$
we see that we have a further factorization 
$$
{\rm DD}^2:\BB^{2n}U(1)_{\nabla}\to \BB^{4n+1}\ZZ/2\to \flat \BB^{4n+1}U(1) \to \BB^{4n+1}U(1)_{\nabla}\;.
$$
Since the DB-cup product refines the classical cup product we can extend this map to a homotopy commutative diagram
$$
\xymatrix{
\BB^{2n}U(1)_{\nabla}\ar[d]^{\rho_2 I}\ar[rr]&&  \BB^{4n+1}\ZZ/2\ar[rr]\ar[drrrr]_{\beta_2} && 
\flat \BB^{4n+1}U(1)\ar[rr] \ar[drr]^{\rho_2 \widetilde{\beta}} && \BB^{4n+1}U(1)_{\nabla}\ar[d]
\\
 \BB^{2n+1}\ZZ/2\ar[rrrrrr]^{\hspace{-1cm}{\iota^2}} &&&&&&\BB^{4n+2}\ZZ/2\;,
 }
 $$
 where, recall, $\widetilde{\beta}$ is the Beckstein corresponding to the exponential sequence. 
By Proposition \ref{Prop mod p} (or Theorem \ref{characterization theorem}) there is an operation 
$\varphi:\BB^{2n}\ZZ/2\to \BB^{4n+1}\ZZ/2$ which fills the left corner
$$
\xymatrix{
\BB^{2n}U(1)_{\rm conn}\ar[d]^{\rho_2 I}\ar[rr]&&  \BB^{4n+1}\ZZ/2\ar[rr]\ar[drrrr]_{\beta_2} &&
 \flat \BB^{4n+1}U(1)\ar[rr] \ar[drr]^{\rho_2 \widetilde{\beta}} && \BB^{4n+1}U(1)_{\rm conn}\ar[d]
\\
 \BB^{2n+1}\ZZ/2\ar@{-->}[urr]^{\varphi}\ar[rrrrrr]^{\hspace{-1cm}{\iota^2}} &&&&&&\BB^{4n+2}\ZZ/2\;,
 }
 $$
such that everything commutes up to homotopy. The homotopy commutativity of the bottom triangle, along with the fact that $Sq^{2n+1}(\iota)=\iota^2$, implies that $\varphi=Sq^{2n}$. The homotopy commutativity of the top part of the diagram proves the claim.
 \endofproof

 \subsection{Properties of the differential Steenrod operations} 
\label{Sec Prop}

 We now discuss general properties of the differential Steenrod squares. These properties can be directly  
  deduced from the general form of these operations as
 \(
 \widehat{Sq}^{2m+1}=j\Gamma_2 Sq^{2m} \rho_2 I\;,
\label{Sq hat}
 \)
 but we make them explicit for the sake of completeness.

\begin{theorem}
[Properties of differential Steenrod squares] The operations $\widehat{Sq}$ satisfy the following:

\begin{enumerate}
\item {\rm Refinement:} The mod 2 reduction of the integral class $I\widehat{Sq}^{2m+1}$ is $Sq^{2m+1}\rho_2I$. 
\item {\rm Torsion}: $\widehat{Sq}^{2m+1}$ takes values in 2-torsion.
\item {\rm  Connectivity:} $\widehat{Sq}^{2m+1}=0, q<0$.
\item {\rm Linearity:} $\widehat{Sq}^{2m+1}(\hat{x}+\hat{y})=\widehat{Sq}^{2m+1}(\hat{x})+\widehat{Sq}^{2m+1}(\hat{y})$.
\item {\rm Squaring:} For $\hat{x}$ of degree $2n+1$, we have $\widehat{Sq}^{2n+1}(\hat{x})=\hat{x}^2$.
\item {\rm Finiteness:} $\widehat{Sq}^{2m+1}(\hat{x})=0$ when $q>n$.
\item {\rm Adem relations:} For even integers $a$ and $b$, we have 
$$
\widehat{Sq}^a\widehat{Sq}^b= \sum_{c\ {\rm odd}} \binom{b-c-1}{a-2c}\widehat{Sq}^{a+b-c}\widehat{Sq}^c\;.
$$
\end{enumerate}
\label{Thm Adem}
\end{theorem}

\theproof
We have already proven properties (1), (2), (3), (5) and (6). For property (4), we have
$$
\widehat{Sq}^{2m+1}(\hat{x}+\hat{y}) = j\Gamma_2 Sq^{2m} \rho_2 I (\hat{x}+\hat{y})\;.
$$
Since the morphisms $I$, $\rho_2$, $\Gamma_2$ and $j$ are induced from homomorphisms of abelian 
groups, they represent linear operations. Since the classical Steenrod squares are linear, we
 immediately deduce that the right hand side is 
$$
j\Gamma_2 Sq^{2m} \rho_2 I (\hat{x}+\hat{y})=j\Gamma_2 Sq^{2m} \rho_2 I (\hat{x})+j\Gamma_2 Sq^{2m} \rho_2 I (\hat{y})\;,
$$
which gives the result. To prove property (7), we have that, by definition,
$\widehat{Sq}^a=j\Gamma_2 Sq^{a-1} \rho_2 I$.
Since $j\Gamma_2$ is linear and takes values in 2-torsion, we have
\bea 
j\Gamma_2 \Big( \sum_c \binom{b-c-1}{a-2c-1}_2
 Sq^{a+b-c-1}Sq^{c} \Big)\rho_2I 
&=& \sum_c \binom{b-c-1}{a-2c-1}_2 \rho_2Ij\Gamma_2 Sq^{a+b-c-1}Sq^{c} \rho_2 I
\\
&=& \sum_{c \ {\rm odd}} \binom{b-c-1}{a-2c-1}j\Gamma_2 Sq^{a+b-c-1}\rho_2Ij\Gamma_2Sq^{c-1}\rho_2I\;.
\\
&=& \sum_{c \ {\rm odd}} \binom{b-c-1}{a-2c-1}\widehat{Sq}^{a+b-c}\widehat{Sq}^{c}\;.
\eea
Recall that we have the relation on binomial coefficients
$$
\binom{b-c-2}{a-2c}=\binom{b-c-1}{a-2c-1}+\binom{b-c-1}{a-2c}\;.
$$
The left hand side is easily seen to be $0\ {\rm mod}\ 2$. Hence 
$$
\binom{b-c-1}{a-2c-1}=\binom{b-c-1}{a-2c}\ {\rm mod}\ 2.
$$ 

\vspace{-6mm}
\endofproof

\begin{remark} [No identity] 
Note that there is no refined Steenrod square which acts 
as an identity. One might be tempted to say $\widehat{Sq}^1$ is 
 ${\rm Id}$ by looking at the definition (see expression \eqref{Sq hat}) 
 and the fact that 
 the identity on the classical Steenrod square is $Sq^0$. However, the effect of 
other morphisms acting on this `classical identity' from both sides 
will affect it nontrivially. Indeed, from the diagram in Proposition 
\ref{Prop trapezoid}, we have the following effect on a differential class
$\hat{x}$
$$
\xymatrix{
\hat{x} \ar @{|->}[r]  & \rho_2 I (\hat{x})  \ar @{|->}[r]^{Sq^0} & 
\rho_2 I (\hat{x}) \ar @{|->}[rrr]^{\beta_2=\rho_2 I \circ j \circ \Gamma_2} &&& 
\beta_2 \rho_2 I (\hat{x}) 
}\;,
$$
which is equal to $Sq^1 \rho_2 I (\hat{x})$, evidently not the identity. 
 Of course this is just a confirmation 
that $\widehat{Sq}^1$ is  a refinement of $Sq^1$.  
This lack of identity is one of several reasons why there should be no notion of 
a differential Steenrod algebra. One could also make such a statement already at
the level of integral cohomology where the identity is certainly not 2-torsion. 
\end{remark}

\begin{remark} [No Cartan formula]
A Cartan formula, which would expand $\widehat{Sq}(\hat{x} \cup_{DB} \hat{y})$
in terms of the product of $\widehat{Sq}(\hat{x})$
 and $\widehat{Sq}(\hat{x})$, with appropriate combinations of degrees, 
 does not seem to exist. 
 There are several reasons for this. The most immediate is that the 
 even Steenrod squares admit no refinements (Proposition \ref{Prop even}). 
 Explicitly, in expanding an 
 odd Steenrod square into a sum of products, each one of the 
 summands will be necessarily a product that  involves an 
 even differential Steenrod square, which does not exist. That is,
 it boils down to the basic fact that an odd number (the degree of the Steenrod
 square) cannot be split into a sum of two odd numbers. 
 \end{remark}

\begin{remark} [Steenrod reduced powers]
Similarly, for odd primes one has an analogue  of Theorem \ref{Thm Adem} for the classical
 Steenrod reduced powers. We do not spell this out, as the 
proofs will be very similar to  those of the above theorem, with 
 obvious changes to coefficients. 
\end{remark}

\subsection{Applications}
\label{Sec App}

In this section, we offer several applications. It will be particularly useful to describe the differential 
Steenrod squares in terms of bundle data as done in \cite{GS1} for the  refined Massey products. 
In order to do this, we will need to make use of the stability of differential Steenrod squares under 
delooping. Note that differential cohomology does not obey the strict suspension isomorphism, i.e. 
in general
$$
\widehat{H}^n(X;\ZZ)\not\simeq \widehat{H}^{n+1}(\Sigma X;\ZZ)\;.
$$
However, we will see that there is a particular sense of stability which 
the refined Steenrod operations enjoy. Let us consider this in more precise terms. 
For an odd integer $2n+1$, the squaring operation is equal to the top Steenrod square
\(
\label{deloop}
\widehat{Sq}^{2n+1}=\cup^2:\BB^{2n}U(1)_{\nabla}\to \flat\BB^{4n+1}U(1)_{\nabla}\;.
\)
Now we can deloop this map $k$ times to get an operation
$$
\BB^k\big(\cup^2\big):\BB^k\big(\BB^{2n}U(1)_{\nabla}\big)\to 
\BB^{k}\big(\flat\BB^{4n+1}U(1)_{\nabla}) \simeq \flat\BB^{4n+1+k}U(1)_{\nabla}\;.
$$
Note that delooping does not commute with having a connection in general, and this is 
indicated by  the parentheses.
\footnote{This view on stacks has been used in \cite{FRS}, where interesting results on 
geometric quantization arise by considering 
e.g. $\BB(\BB U(1)_{\rm conn})$ as the 2-stack whose sections are $U(1)$-bundle 
gerbes with connective structure but without curving, and in \cite{Al} to 
succinctly get results on String structures that otherwise require considerable buildup.} 
However, the two operations commute when 
the stack is flat, as is indicated in the equivalence on the right hand side.
Explicitly, we can write down what this map does on sections (at the level of 
chain complexes) using the formula for the DB-cup product; the map 
$\BB^k\big(\cup^2\big)$ is given by
$$
\BB^k\big(\cup^2\big)(\alpha)=\left\{\begin{array}{ccl}
\alpha^2 && \text{if}\ {\rm deg}(\alpha)=2n+k+1,
\\
\alpha\wedge d\alpha && \text{if}\  {\rm deg}(\alpha)=2n,
\\
0 && \text{otherwise}.
\end{array}
\right.
$$
Now the integration map $I$ and mod 2 reduction $\rho_2$ clearly commute with delooping. 
Furthermore, since the differential Steenrod squares refine the classical one, we have
$$
\rho_2I\BB^k\big(\cup^2\big)=\rho_2I\BB^k\big(\widehat{Sq}^{2n+1}\big)=\BB^k\big(\rho_2I\widehat{Sq}^{2n+1}\big)
=\BB^k\big(Sq^{2n+1}\rho_2I\big)\;.
$$ 
Then stability of the classical Steenrod squares implies that the right hand side 
is just $Sq^{2n+1}\rho_2I$. Hence, the stacky delooping $\BB^k\big(\widehat{Sq}^{2n+1}\big)$ refines the classical operation $Sq^{2n+1}$. However, the source of this operation is not the stack $\BB^{2n+k}U(1)_{\nabla}$. To remedy this, we simply note that we have a (in fact, strictly) commutative diagram
$$
\xymatrix{
\BB^{2n+k}U(1)_{\nabla} \ar[rr]^{\rm pr}\ar[dr]_{I} && \BB^k\big(\BB^{2n}U(1)_{\nabla}\big)
\ar[dl]^{I}
\\
& \BB^{2n+k+1}\ZZ &
}\;,
$$
where ${\rm pr}$ is the projection map induced from the morphism of chain complexes
$$
\xymatrix{
\ZZ \ar[d] \ar[r] & \Omega^0 \ar[r] \ar[d] & \cdots \ar[r] & \Omega^{2n} \ar[r]  \ar[d] & 
\Omega^{2n+1} \ar[r]  \ar[d] & 
\cdots \ar[r] & \Omega^{2n+k} \ar[d]
\\
\ZZ \ar[r] & \Omega^0 \ar[r] & \cdots \ar[r] & \Omega^{2n} \ar[r]  & 0 \ar[r]  & 
\cdots \ar[r] & 0\;.
}
$$
From the commutativity, we find the following.
\begin{proposition} 
\label{Prop stable} 
The refined Steenrod squares are stable, in the sense that they give morphisms of stacks
$$
\widehat{Sq}^{2n+1}=\BB^k\big(\widehat{Sq}^{2n+1}\big){\rm pr}: ~\BB^{2n+k}U(1)_{\nabla}\to \flat \BB^{4n+k+1}U(1)_{\nabla}\;.
$$
\end{proposition} 

\medskip
The next example will present the general case. Then we will narrow our scope to specific cases arising from physics,
which can be seen as a continuation and extension of the discussions in 
 \cite{FSS1} \cite{FSS2} \cite{E8} 
\cite{GS1}.

\begin{example}
Let   
$$
\hat{x}:X\to \BB^{2n+k}U(1)_{\nabla}
$$
be a $(2n+k-1)$-bundle representing a differential cohomology class and let $\{U_{\alpha}\}$ be a good open cover of $X$. 
Then, from \cite{FSSt}, $\hat{x}$ determines the following data:
\begin{itemize}
\item An integral {\v C}ech cocycle $n_{\alpha_0\hdots \alpha_{2n+k+1}}$ on $(2n+k+1)$-fold intersections.
\item Differential $j$-forms, $j\leq 2n+k$, ${\cal A}_{j}$ on $(2n+k-j)$-fold intersections such that
$$(-1)^k\delta {\cal A}_{j}=d{\cal A}_{j-1}\;.$$
\end{itemize}
Post-composing the map $\hat{x}$ with the map \eqref{deloop}, gives the refined degree $(2n+1)$ Steenrod square
$$
\widehat{Sq}^{2n+1}(\hat{x}):X \to \BB^{2n+k}U(1)_{\nabla}\to \BB^k\big(\BB^{2n}U(1)_{\nabla}\big)
\to \flat\BB^{4n+k+1}U(1)_{\nabla}\;.
$$
In terms of bundle data, this composite is fairly straightforward to write down explicitly. The combinatorics involved can be deduced from a rather tedious (but straightforward) calculation involving the analogue of the map in unbounded sheaves of chain complexes. Explicitly, we get the following bundle data.

\begin{itemize}
\item An integral {\v C}ech cocycle 
$$n_{\alpha_0\hdots \alpha_{2n+k}}n_{\alpha_{2n}\alpha_{2n+1}\hdots \alpha_{4n+k}}$$
on $(4n+k+1)$-fold intersections.
\item Differential $j$-forms, $0\leq j\leq 2n$, $n_{\alpha_0\hdots \alpha_{2n+k}}{\cal A}_j$ on $(4n+k+1-j)$-fold intersections satisfying the {\v C}ech-Deligne cocycle condition.
\item Differential $(j+2n+1)$-forms, $0\leq j\leq 2n$, ${\cal A}_j\wedge d{\cal A}_{2n}$ on $(2n+k-j)$-fold intersections satisfying the {\v C}ech-Deligne cocycle condition.
\item $0$ on all other intersections.
\end{itemize}
\end{example}

\medskip
\noindent The following two examples are directed related, but we split for ease of presentation. 
We will calculate $\widehat{Sq}^1$ and $\widehat{Sq}^3$ in terms of bundle data.

\begin{example} 
\label{Ex H}
Let $H$ be a closed 3-form. Suppose that $H$ admits a differential refinement $\widehat{H}$. Then we can identify $\widehat{H}$ with a bundle given by the following data:
\begin{itemize}
\item An integral {\v C}ech cocycle $n_{\alpha\beta\gamma\delta}$ on quadruple intersections.
\item A smooth, real valued function $f_{\alpha\beta\gamma}$ on triple intersections (satisfying the {\v C}ech-Deligne cocycle condition).
\item A differential 1-form ${\cal A}_{\alpha\beta}$ on intersections (satisfying the {\v C}ech-Deligne cocycle condition).
\item A differential 2-form ${\cal B}_{\alpha}$ on open sets (satisfying the {\v C}ech-Deligne cocycle condition).
\end{itemize}
The refined Steenrod square $\widehat{Sq}^3$, being in top degree, is the usual Deligne-Beilinson square given by the data 
\begin{itemize}
\item An integral {\v C}ech cocycle $n_{\alpha\beta\gamma\delta}n_{\delta\epsilon\eta\xi}$ on 7-fold intersections.
\item A smooth, real valued function $n_{\alpha\beta\gamma\delta}f_{\delta\epsilon\eta}$ on 6-fold intersections.
\item A differential 1-form $n_{\alpha\beta\gamma\delta}{\cal A}_{\delta\epsilon}$ on 5-fold intersections.
\item A differential 2-form $n_{\alpha\beta\gamma\delta}{\cal B}_{\delta}$ on 4-fold intersections.
\item $f_{\alpha\beta\gamma}\wedge H$ on 3-fold intersections.
\item ${\cal A}_{\alpha\beta}\wedge H$ on intersections.
\item ${\cal B}_{\alpha}\wedge H$ on open sets.
\end{itemize}
\end{example}

\begin{example}
Continuing Example \ref{Ex H},  
$\widehat{Sq}^1(\widehat{H})$ is given by
\begin{itemize}
\item An integral {\v C}ech cocycle $n_{\alpha\beta\gamma\delta}n_{\beta\gamma\delta \xi}$ on 5-fold intersections.
\item A smooth, real valued function $n_{\alpha\beta\gamma\delta}f_{\beta\gamma\delta}$ on 4-fold intersections.
\item $f_{\alpha\beta\gamma}\wedge df_{\alpha\beta\gamma}$ on 3-fold intersections.
\item $0$ on all other intersections.
\end{itemize}
\end{example}

\begin{remark}[Global patterns]
{\bf (i)} In each of the bundles in the last two examples, it is interesting to note the form of the
 first nonzero value. These are
\(
f_{\alpha\beta\gamma}\wedge df_{\gamma\epsilon\eta},\qquad
 {\cal B}_{\alpha}\wedge H={\cal B}_{\alpha}\wedge dB_{\alpha}\;,
\label{CS}
\)
respectively. From left to right these correspond to the first,
second, and third refined Steenrod operation, respectively. 

\vspace{1mm}
\noindent  {\bf (ii)} In fact, an interesting pattern emerges when one 
looks at the first non-zero information in the Steenrod squares. This is summarized in the following table.

\medskip
\begin{center}
\begin{tabular}{| l | c | c | c | c |}
\hline
Degree of curvature form & ${\rm deg}(H)=1$ & ${\rm deg}(H)=3$ & ${\rm deg}(H)=5$ & ${\rm deg}(H)=7$
\\ \hline
Degree of square of $H$ & 2 & 6 & 10 & 14
\\ \hline
Degree of first nonzero form in $\widehat{Sq}^1(\widehat{H})$ & 1-form  & 1-form & 1-form & 1-form 
\\
Dimension of bundle  & 1 & 3 & 5 & 7  
\\ \hline
Degree of first nonzero form in $\widehat{Sq}^3(\widehat{H})$ & N/A & 5-form & 5-form & 5-form 
\\
Dimension of bundle & & 5 & 7 & 9  
\\ \hline
Degree of first nonzero form in $\widehat{Sq}^5(\widehat{H})$ & N/A & N/A & 9-form & 9-form 
\\
Dimension of bundle  & & & 9 &11  
\\ \hline
Degree of first nonzero form in $\widehat{Sq}^7(\widehat{H})$ & N/A & N/A & N/A & 13-form 
\\
Dimension of bundle  & & & & 13 
\\ \hline

\end{tabular}
\end{center}
\end{remark}

\medskip
\begin{remark}[Physical action functionals]
The proper way to formulate action functionals is to to do so with a view
towards towards quantization. That is, the exponentiated action 
functional, which is the integrand in a path integral over the configuration/moduli
space needed for quantization, should be well-defined. The above expressions \eqref{CS},
are of Chern-Simons type and automatically take values in $\R/\Z$ by our constructions,
hence give rise to well-defined path integrands. These are also secondary 
classes associated with flat bundles. The discussions in \cite{GS1}, being about
Massey products in differential cohomology, can then be viewed as 
a secondary formulation of a secondary formulation. The action functionals
reflect that effect.  In both formulations the main point, as far as action
functionals go, is that the cup product being zero does not end the story, but
rather opens up the possibility for a considerable amount of geometry and dynamics 
to be captured by resorting to secondary considerations. 
 \end{remark}

\medskip 
We close by noting that there are many possible further applications to geometry and topology. 
Some of these discussions will appear soon. Indeed, the first application will
be to the Atiyah-Hirzebruch spectral sequence in differential generalized cohomology theories
\cite{GS3}.

\vspace{6mm}
\noindent {\bf \large Acknowledgement}

\medskip
\noindent The authors would like to thank Ulrich Bunke and Thomas Nikolaus
for interesting discussions at the early stages of this project.  The authors would also like to thank Thomas Schick for pointing out an error in the K\"unneth formula in the published version of this article.


\end{document}